\documentclass[10pt, a4paper]{article}

\usepackage{vmargin}

\setpapersize{A4}
\setmargins{2.5cm}       
{1.5cm}                        
{16.5cm}                      
{23.42cm}                    
{10pt}                           
{1cm}                           
{0pt}                             
{2cm}                           
\date{}

\usepackage{enumerate} 

\usepackage{amsmath,amssymb,latexsym,
}
\usepackage{latexsym,amssymb}
\usepackage{amsmath}
\usepackage{mathrsfs}
\usepackage{bussproofs}
\usepackage{float}
\usepackage{tikz}
\usepackage{amsfonts}   

\newtheorem{defi}{Definition}[section]
\newtheorem{theo}[defi]{Theorem}
\newtheorem{lem}[defi]{Lemma}

\newtheorem{prop}[defi]{Proposition}
\newtheorem{rem}[defi]{Remark}

\newtheorem{cor}[defi]{Corollary}

\def\wff{\mathop{\rm WFF}\nolimits}
\def\sent{\mathop{\rm SENT}\nolimits}
\newcommand{\negr}[1]{\boldsymbol{#1}}

\newenvironment{dem}{\noindent\bf Proof. \rm}{\hfill $\negr{\blacksquare}$}

\newcommand{\pt}{\forall} 
\newcommand{\ex}{\exists} 
\newcommand{\no}{\neg} 
\newcommand{\cons}{{\circ}} 

\def\imp{\!\to\!}
\def\sii{\!\leftrightarrow\!}

\def\fun{\! \longrightarrow\!}
\def\d{\displaystyle}

\newcommand{\mm}{\mathcal} 
\newcommand{\rr}{\mathscr} 
\newcommand{\bb}{\mathbb}  
\newcommand{\ff}{\mathfrak} 

\newcommand{\tripc}[1]{\|#1\|^{\mathfrak A}}
\newcommand{\ese}[1]{S(\ff #1)}

\newcommand{\mas}[1]{\|#1\|^{\mathfrak A}_\oplus}
\newcommand{\por}[1]{\|#1\|^{\mathfrak A}_\odot}
\newcommand{\menos}[1]{\|#1\|^{\mathfrak A}_\ominus}

\newcommand{\ciore}{{\bf Ciore}}
\newcommand{\qciore}{{\bf QCiore}}
\newcommand{\quci}{{\bf QCiore}}

\newcommand{\gciore}{{\bf GCiore}}
\newcommand{\gqciore}{{\bf GQCiore}}
\newcommand{\lfis}{{\bf LFI}s}

\newcommand{\lfidos}{{\bf LFI2}}

\newcommand{\brac}[1]{\langle #1\rangle}	
\newcounter{defcounter}
\newenvironment{myequationt}{%

\begin{equation}}
{\end{equation}}
\newenvironment{myequationb}{%

\begin{equation}}
{\end{equation}}

\newenvironment{myequationk}{%

\begin{equation}}
{\end{equation}}

\newenvironment{myequationbc1}{%

\begin{equation}}
{\end{equation}}

\newenvironment{myequationci}{%

\begin{equation}}
{\end{equation}}

\newenvironment{myequationcf}{%

\begin{equation}}
{\end{equation}}

\newenvironment{myequationco1}{%

\begin{equation}}
{\end{equation}}
\newenvironment{myequationco2}{%

\begin{equation}}
{\end{equation}}
\newenvironment{myequationco3}{%

\begin{equation}}
{\end{equation}}

\newenvironment{myequation}{%
\addtocounter{equation}{-1}
\refstepcounter{defcounter}

\begin{equation}}
{\end{equation}}

\newcommand{\ter}[1]{t_{ #1}^\ff A[s]}

\def\um{\mathbf{\frac 1 2}}
\def\0{\mathbf{0}}
\def\1{\mathbf{1}}

\def\A{\mathcal{A}}
\def\B{\mathcal{B}}

\def\T{\mathbf{T}}
\def\C{\mathcal{C}}
\def\D{\mathcal{D}}

\def\red{\color{red}}

\title{Proof-theoretic aspects of paraconsistency with strong consistency operator}

\author{Victoria Arce Pistone\footnote{The first author is grateful to CONICET for providing financial support for this research.} and Mart\'in Figallo \footnote{The second author was partially supported by the {\em Visiting Researcher Award} program funded by FAPESP grant 2022/03862-2.}}

\date{\textit{\small Departamento de Matem\'atica  and Instituto de Matem\'atica (INMABB). Universidad Nacional del Sur. Bah\'{\i}a Blanca, Argentina}}

\begin{document}

\maketitle

\begin{abstract}
 In order to develop efficient tools for automated reasoning  with inconsistency (theorem provers), eventually making Logics of Formal inconsistency ({\bf LFI}) a more appealing formalism for reasoning under uncertainty, it is important to develop the proof theory of the  first-order versions of such {\bf LFI}'s.  In our work, we intend make a first step in that direction. On the other hand, the logic \ciore\ was developed to provide new logical systems in the study of inconsistent databases from the point of view of {\bf LFI}. An interesting fact about \ciore\ is that it has a {\em strong} consistency operator, that is, a consistency operator which (forward/backward) propagates inconsistency. Also, it turns out to be an algebraizable logic (in the sense of Blok and Pigozzi) that can be characterized by means of a 3-valued logical matrix.  Recently, a first-order version of \ciore, namely \qciore, was defined  preserving the spirit of \ciore, that is, without introducing unexpected relations between the quantifiers. Besides, some important model-theoretic results were obtained for this logic.

In this paper we study some proof--theoretic aspects of both \ciore\ and \qciore\ respectively. In first place, we introduce a two-sided sequent system for \ciore. Later, we prove that this system enjoys the cut-elimination property and apply it to derive some interesting properties. Later, we extend the above-mentioned system to first-order languages and prove completeness and cut-elimination property using the well-known Sh\"utte's technique. 

\vspace*{4mm}


\vspace*{4mm}


\end{abstract}

\section{Introduction}
A paraconsistent logic is a formal system that allows to reason about inconsistent information without lapsing into absurdity. In a non-paraconsistent setting, inconsistency {\em explodes} in the sense that if a contradiction obtains,  then everything obtains. The first systematic study of paraconsistent logics was carried out by da Costa, when presented in~\cite{daC} his well-known hierarchy $C_n$ (for $n \geq 1$) of systems. His approach to paraconsistency, nowadays known as the {\em Brazilian school of paraconsistency}, was naturally generalized by W. Carnielli and J. Marcos in~\cite{Tax} with the notion of  {\em Logics of Formal Inconsistency} (\lfis, for short).  These are paraconsistent logics that internalize the very notions of consistency and inconsistency at the object-language level. In~\cite{Tax}, an important subclass of \lfis\ , called {\bf C}-systems, was considered. Similarly to \lfis\ ,  {\bf C}-systems are built over the positive basis of some given consistent logic and with a special connective $\circ$ (either primitive or defined) that allows to express the notion of consistency of sentences inside the object  language.

As it was pointed out in \cite{Av00}, since their introduction in terms of Hilbert-style systems in the 1960s, the main obstacle to efficient use of {\bf C}-systems has been the lack of analytic calculi for them. In this same paper, the authors provide a uniform and modular method for a systematic generation of cut-free sequent calculi for
a large family of paraconsistent logics. All the systems studied by them have semantics in terms of non-deterministic matrices (Nmatrices) which are a natural generalization of standard multi-valued matrices obtained by importing the notion of non-deterministic computations from computer science into the truth-tables of logical connectives. They use a method from \cite{Avron02} for constructing cut-free Gentzen-type systems for logics which have a characteristic finite-valued non-deterministic matrices (Nmatrices) and whose language is sufficiently expressive, in a certain sense.

On the other hand, the 3-valued paraconsistent logic \ciore\ was developed by Carnielli, Marcos and de Amo under the name \lfidos\ in the study of inconsistent databases from the point of view of {\em Logics of Formal Inconsistency} ({\bf LFI}s). They studied this logic considering a primitive inconsistency connective $\bullet$ instead of a consistency connective $\circ$. As observed in \cite{Tax}, contrary to $C_1$ and $C^+_1$,  \ciore\ is algebraizable in the sense of Blok and Pigozzi. More than this, because of the strong properties enjoyed by the consistency connective $\circ$ and by the paraconsistent negation, it can be characterized by means of a 3-valued logical matrix. 

It is important to mention that \ciore\ does not belong to the family of $C$-systems studied in \cite{Av00}. Indeed, the paraconsistent negation of \ciore\ does not validate the De Morgan laws but it does validate the following very particular law (see \cite[Thm 2.2 (xii)]{ConGomFil})
$$\left((\alpha\wedge\neg\alpha)\wedge(\beta\wedge\neg\beta)\right)\to \neg(\alpha\vee\beta)$$
No $C$-system studied in \cite{Av00} consider this alternative; nor the extreme laws of propagation and retropropagation verified by the operator $\circ$ (strong consistency operator). Carnielli and Marcos also introduced a first-order version of \ciore\ called \lfidos$^*$. As it was noted in \cite{ConGomFil}, \lfidos$^*$ satisfies a somewhat counter-intuitive property: the universal and the existential quantifier are inter-definable by means of the paraconsistent negation, as it happens in classical first-order logic with respect to the classical negation. This feature seems to be unnatural, given that both quantifiers have the classical meaning in \lfidos$^*$, and that this logic does not satisfy the De Morgan laws with respect to its paraconsistent negation. Due to this, in \cite{ConGomFil} it was introduced a first-order version of \ciore ,  named \qciore ,  preserving the spirit of \ciore, that is, without introducing unexpected relationships between the quantifiers. Some important results of classical Model Theory are obtained for this logic, such as Robinson’s joint consistency theorem, amalgamation and interpolation. The main purpose of \cite{Av00} was to develop efficient tools for automated reasoning with inconsistency, eventually making {\bf LFI}'s a more appealing formalism for reasoning under uncertainty. However, it is clear that for the purposes of building {\bf LFI}-based theorem provers for real-life applications, it is important to develop the proof theory of the  first-order versions of such {\bf LFI}'s. In our work, we intend make a first step in that direction. 

\

In this work we study some proof--theoretic aspects of both \ciore\ and \qciore\ respectively. In Section \ref{s2} we recall all notions an results known concerning \ciore\ and \qciore\ as well as we state notation. In Section \ref{s3}, inspired by the method depicted in \cite{Avron02}, we present a sequent-style systems \ciore\ and, since we do not make explicit the different stages of this long process, in Section \ref{s4} we provide a semantical proof of the cut-elimination property and, in Section \ref{s5}, we show some applications. In Section \ref{s6} we extend the above-mentioned system to first-order languages and, using the well-known Sh\"utte's technique, we prove the completeness and cut-elimination theorems. Finally, in Section \ref{s7} we draw some conclusions and describe some future work.

\section{Preliminaries}\label{s2}
Let $\mathscr{L}$ be a propositional language and let $\mathfrak{Fm}$ be the absolutely free algebra over $\mathscr{L}$ generated by some denumerable set of propositional variables, with underlying set (of formulas) $Fm$ and let ${\cal M}=\langle {\cal T}, {\cal D }, {\cal O}\rangle$ be a logic matrix for $\mathscr{L}$, that is, ${\cal T}$ is a finite, non-empty set of truth values, ${\cal D}$ is a non-empty proper set of ${\cal T}$, and ${\cal O}$ includes a $k$-ary function $\hat{f}: {\cal T}^k\to {\cal T}$ for each $k$-ary connective $f$.
Recall that, a valuation in ${\cal M}$ is a function $v:Fm\to {\cal T}$ such that
$$v(f(\psi_1, \dots, \psi_k))=\hat{f}(v(\psi_1), \dots, v(\psi_k))$$
for each $k$-ary connective $f$ and all $\psi_1, \dots, \psi_k\in Fm$. 
A formula $\alpha \in Fm$ is satisfied by a given valuation $v$, in symbols $v\models \alpha$, if $v(\alpha)\in {\cal D}$. A sequent $\Gamma \Rightarrow \Delta$ is satisfied by the valuation $v$, in symbols $v\models \, \Gamma \Rightarrow \Delta$, if either $v$ does not satisfy some formula in $\Gamma$ or $v$ satisfies some formula in $\Delta$. A sequent is {\em valid} (w.r.t the matrix ${\cal M}$) if it is satisfied by all valuations. We write $\vdash_{\cal M}  \Gamma\Rightarrow \Delta$ to indicate that the sequent  $\Gamma\Rightarrow \Delta$ is valid in $\cal M$.\\ 
Now, suppose that ${\cal T}=\{t_0,\dots, t_{n-1}\}$, where $n\geq 2$,  and ${\cal D}=\{t_d,\dots, t_{n-1}\}$, where $1\leq d \leq n-1$. An $n$--sequent (\cite{Avron01})  over $\mathscr{L}$ is an expression  
$$ \Gamma_0 \mid \dots \mid \Gamma_{n-1}$$
where, for each $i$, $\Gamma_i$ is a finite set of formulas. A valuation $v$ satisfies the $n$--sequent $\Gamma_0 \mid \dots \mid \Gamma_{n-1}$ iff there exists $i$, $0\leq i \leq n-1$ \, and $\psi\in \Gamma_i$ such that $v(\psi)=t_i$. An $n$--sequent is valid if it is satisfied by every valuation $v$. It is clear that a valuation $v$ satisfies an ordinary sequent $\Gamma  \Rightarrow \Delta$ iff $v$ satisfies the $n$--sequent $\Gamma_1 \mid \dots \mid \Gamma_{n-1}$ where $\Gamma_i=\Gamma$ for all $0\leq i\leq d-1$ and $\Gamma_j=\Delta$ for all $d\leq j\leq n-1$ . \\
An alternative presentation of $n$-sequents is by means of sets of {\em signed formulas}. A signed formula over the language $\mathscr{L}$ and ${\cal T}$, is an expression of the form
$$t_i : \psi$$
where $t_i\in  {\cal T}$ and $\psi \in Fm$. A valuation $v$ satisfies the signed formula $t_i : \psi$ iff $v(\psi)=t_i$. If $\Omega\subseteq {\cal T}$ and $\Gamma \subseteq Fm$, we denote by $\Omega : \Gamma$ the set
$$\Omega : \Gamma =\{ t:\alpha \mid t\in \Omega, \alpha \in \Gamma\}$$
If $\Omega=\{t\}$, we write $t : \Gamma$ instead of $\{t\} : \Gamma$.  
A valuation satisfies the set of signed formulas $\Omega : \Gamma$  if it satisfies some signed formula of $\Omega : \Gamma$; and we say that $\Omega : \Gamma$ is valid if it is satisfied by every valuation $v\in {\cal V}$.
It is clear that,  the $n$--sequent $ \Gamma_0 \mid \dots \mid \Gamma_{n-1}$ is valid iff the set of signed formulas $\bigcup \limits_{i=0}^ {n-1} t_i : \Gamma_i$ is valid.

\

Let ${\cal L}$ be the propositional language defined over the propositional signature $\Sigma=\{\wedge, \vee, \rightarrow, \neg, \circ\}$. In \cite{ConGomFil}, the propositional logic \ciore\  was presented, over the language ${\cal L}$, by means of the following Hilbert-style system (as usual, $\alpha\leftrightarrow\beta$ denotes the formula $(\alpha\rightarrow\beta)\wedge(\beta\rightarrow\alpha)$, $\alpha, \beta \in Fm$):\\[2mm]
\noindent{\bf Axiom schemata:} 

\vspace{-0.3cm}
\begin{myequation}
  \alpha\imp(\beta\imp\alpha)
   \end{myequation}
\vspace{-0.5cm}
\begin{myequation}
  \big(\alpha\imp(\beta\imp\gamma)\big)\imp\big((\alpha\imp\beta)\imp(\alpha\imp\gamma)\big)
   \end{myequation}
\vspace{-0.5cm}
\begin{myequation}
  \alpha\imp\big(\beta\imp(\alpha\wedge\beta)\big)
   \end{myequation}
\vspace{-0.5cm}
\begin{myequation}
  (\alpha\wedge\beta)\imp\alpha
   \end{myequation}
\vspace{-0.5cm}
\begin{myequation}
  (\alpha\wedge\beta)\imp\beta
   \end{myequation}
\vspace{-0.5cm}
\begin{myequation}
  \alpha\imp(\alpha\vee\beta)
   \end{myequation}
\vspace{-0.5cm}
\begin{myequation}
  \beta\imp(\alpha\vee\beta)
   \end{myequation}
\vspace{-0.5cm}
\begin{myequation}
  \big(
\alpha\imp\gamma
\big)
\imp
\big(
(\beta\imp\gamma)\imp((\alpha\vee\beta)\imp\gamma)\big
)
   \end{myequation}
\vspace{-0.5cm}
\begin{myequation}
  (\alpha\imp\beta)\vee\alpha
   \end{myequation}
\vspace{-0.5cm}
\begin{myequation}
  \alpha\vee\neg\alpha
   \end{myequation}
\vspace{-0.5cm}
\begin{myequationbc1}
\circ\alpha\imp\big(\alpha\imp(\neg\alpha\imp\beta)\big)
\end{myequationbc1}
\vspace{-0.5cm}
\begin{myequationci}
\neg\circ\!\alpha\imp(\alpha\wedge\neg\alpha)
\end{myequationci}
\vspace{-0.5cm}
\begin{myequationcf}
\neg\neg\alpha\sii\alpha
\end{myequationcf}
\vspace{-0.5cm}
\begin{myequationco1}
(\circ\alpha\vee\circ\beta)\sii\circ(\alpha\wedge\beta)
\end{myequationco1}
\vspace{-0.5cm}
\begin{myequationco2}
(\circ\alpha\vee\circ\beta)\sii\circ(\alpha\vee\beta)
\end{myequationco2}
\vspace{-0.5cm}
\begin{myequationco3}
(\circ\alpha\vee\circ\beta)\sii\circ(\alpha\imp\beta)
\end{myequationco3}

\noindent{\bf Inference rule:}
\begin{center}
{\bf (MP)}\hfill $\displaystyle\frac{\alpha \hspace{5mm}\alpha\imp\beta}{\beta}$\hfill\hspace{1mm}
\end{center}
The (Tarskian) consequence relation obtained from the Hilbert calculus for \ciore\  will be denoted by $\vdash_\ciore$. Observe that axioms {\bf (Ax1)-(Ax9)} plus {\bf (MP)} constitute a Hilbert calculus for positive classical logic ${\bf (CPL^{+})}$. As observed above, \ciore\ is algebraizable in the sense of Blok and Pigozzi. Moreover, it can can be characterized by a 3-valued logical matrix.
\begin{theo}\label{SoundCompFI2}
The system \ciore\ is sound and complete with respect to the following three-valued matrix $\mm M_{e}=\langle \mm T,  \mm D, \{\hat\wedge,\hat\vee,\hat\rightarrow, \hat\neg, \hat\circ \}\rangle$ over the signature $\Sigma$ with domain $\mm T=\{\1, \um, \0\}$ and set of
designated values $\mm D=\{\1, \um\}$ such that the truth-tables associated to each connective are the following:

\begin{center}
\begin{tabular}{| c || c | c | c | }
\hline 
$\hat\wedge$ & $\1$ & $\um$ & $\0$ \\ \hline \hline
$\1$ & $\1$ & $\1$ & $\0$\\ \hline
$\um$ & $\1$ & $\um$ & $\0$\\ \hline
$\0$ & $\0$ & $\0$ & $\0$\\ \hline
\end{tabular} \,
\begin{tabular}{| c || c | c | c | }
\hline
$\hat\vee$ & $\1$ & $\um$ & $\0$ \\ \hline \hline
$\1$ & $\1$ & $\1$ & $\1$\\ \hline
$\um$ & $\1$ & $\um$ & $\1$\\ \hline
$\0$ & $\1$ & $\1$ & $\0$\\ \hline
\end{tabular} \, 
\begin{tabular}{| c || c | c | c | }
\hline
$\hat\rightarrow$ & $\1$ & $\um$ & $\0$ \\ \hline \hline
$\1$ & $\1$ & $\1$ & $\0$\\ \hline
$\um$ & $\1$ & $\um$ & $\0$\\ \hline
$\0$ & $\1$ & $\1$ & $\1$\\ \hline
\end{tabular} \,
\begin{tabular}{| c || c || c |}
\hline
 & $\hat\neg$ & $\hat\circ$  \\ \hline\hline
$\1$ & $\0$ & $\1$ \\ \hline
$\um$ & $\um$ & $\0$\\ \hline
$\0$ & $\1$ & $\1$ \\ \hline
\end{tabular}
\end{center}
\end{theo}

In \cite{Av00}, a general method for constructing cut-free sequent calculi for $C$-systems was provided. This  method applies to a large family of $C$-systems, covering many $C$-systems studied in the literature. However, it does not apply to \ciore . Indeed, in \cite{Av00}, the authors consider $C$-systems which are extensions of  {\bf BK} which, in turn, is obtained by adding to the standard Hilbert-style system for the positive fragment (i.e. $\{\vee,\wedge,\to\}$-fragment) of classical propositional logic the axioms
\begin{myequationt}\alpha\vee\neg\alpha
\end{myequationt}
\begin{myequationb}\circ\alpha\to((\alpha\wedge\neg\alpha)\to\beta)
\end{myequationb}
\begin{myequationk}\circ\alpha\vee(\alpha\wedge\neg\alpha)
\end{myequationk}
As we mentioned above, the negation of \ciore\ does not validate the De Morgan Laws. However, it validates the very particular law (see \cite[Thm 2.2 (xii)]{ConGomFil})
$$\left((\alpha\wedge\neg\alpha)\wedge(\beta\wedge\neg\beta)\right)\to \neg(\alpha\vee\beta)$$
No extension studied in \cite{Av00} consider this alternative; nor the extreme laws of propagation and retropropagation verified by the operator $\circ$ which are substantiated by axioms ({\bf cor}$_1$)--({\bf cor}$_3$).

\

In \cite{ConGomFil}, the logic \qciore\ was introduced as a natural first-order version of \ciore\ and with semantics based on the notion of triples (or partial relations). Recall that a first-order signature $\Theta=\brac{\mm P,\mm F,\mm C}$ is composed by: a set $\mm P=\underset{n\in\bb N}{\bigcup} P_n$ such that, for each $n\geq 1$, $\mm P_n$ is a set of predicate symbols of arity $n$; a set $\mm F=\underset{n\in\bb N}{\bigcup} F_n$ such that, for each $n\geq 1$, $\mm F_n$ is a set of function symbols of arity $n$; and a set $\mm C$ of individual constants.

Let $\rr L(\Theta)$ (or just $\rr L$) be the first-order language defined as  usual from the connectives $\wedge$, $\vee$, $\imp$, $\no$, $\cons$, the quantifiers $\pt$, $\ex$,  a denumerable set of free variable symbols ${\mm V}_f=\{a_1, a_2, \dots\}$, a  denumerable set of bound variable symbols: ${\mm V}_b=\{x_1, x_2, \dots\}$ and a given first-order signature $ \Theta$. We denote by $\wff(\rr L(\Theta))$ and $\sent(\rr L(\Theta))$ the set of well-formed formulas and  sentences (formulas without free-variables) over the signature $\Theta$, respectively. 
Any finite sequence of symbols from the language $\rr L$ is an {\em expression} of $\rr L$. Recall that,  if $A$ is an expression and $\tau_1,\dots,\tau_n$ are distinct primitive symbols, and $\sigma_1, \dots, \sigma_n$ are any symbols, then by

$$\left( A \, \d \frac{\tau_1,\dots,\tau_n} {\sigma_1, \dots, \sigma_n} \right)$$

\noindent we mean the expression obtained from $A$ by writing $\sigma_1, \dots, \sigma_n$ in place of $\tau_1,\dots,\tau_n$, respectively, at each occurrence of $\tau_1,\dots,\tau_n$ and where the symbols are replaced simultaneously. Recall, also, that if $A$ is a formula and $t_1, \dots, t_n$ are terms and there is a formula $B$ and free variable symbols $b_1,\dots,b_n$ such that $A$ is $(B\,\frac{b_1,\dots,b_n}{t_1, \dots, t_n})$  we write $B$ as $B(b_1,\dots, b_n)$ and $A$ as $B(t_1, \dots, t_n)$ and we say that for each $i$, $1\leq i \leq n$, the occurrences of the term $t_i$ are {\em indicated} in $A$. The term $t$ is {\em fully indicated}  in $A$ if every occurrence of $t$ is obtained by such replacement (see \cite{Tau}). In what follows, $\Theta$ is a first-order signature.

\begin{defi} (\cite{ConGomFil}). The logic \qciore\ on the language $\rr L(\Theta)$ is defined as the Hilbert calculus obtained by extending \ciore (expressed in the language $\rr L(\Theta)$) by adding the following:

\

\noindent{\bf Axiom schemata:}
\begin{myequation}
  \varphi(t)\imp\ex x\varphi(x), 
   \end{myequation}
\vspace{-0.5cm}

\begin{myequation}
  \pt x\varphi(x)\imp\varphi(t), 
   \end{myequation}
\vspace{-0.5cm}

\begin{myequation}
\cons\ex x\varphi\sii\ex x\cons\varphi,   
\end{myequation}
\vspace{-0.5cm}

\begin{myequation}
\cons\pt x\varphi\sii\ex x\cons\varphi,  
   \end{myequation}

\noindent where $t$ is an arbitrary term; $\varphi(t)$ and $\varphi(x)$ are $\left(\varphi(a)\, \d \frac{a}{t}\right) \mbox{ and } \left(\varphi(a)\, \d \frac{a}{x}\right)$, respectively.

\

\noindent{\bf Inference rules}
\begin{align*}
(\pt\mbox{-}{\bf In}) \ \frac{\phi\rightarrow\psi(a)}{\phi\rightarrow\forall x\psi(x)}
&& \mbox{if $a$ does not occur in $\phi$}\\
(\ex\mbox{-}{\bf In}) \ \frac{\phi(a)\rightarrow\psi}{\exists x \phi(x)\rightarrow\psi}
&& \mbox{if $a$ does not occur in $\psi$}
\end{align*}
 \noindent and where $\varphi(x)$ is $\left(\varphi(a)\, \d \frac{a}{x}\right)$.
\end{defi}

\

In \cite{ConGomFil}, it was introduced a semantic version of \qciore\ by using {\em partial relations} defined in terms of triples. The use of triples proved to be very useful since it allowed to write simpler proofs.

\

Let ${\bf 3}=\{\0,\um, \1\}$ and consider the algebraic structure $\mathsf{3}=\langle {\bf 3}; \vee, \wedge, \imp, \neg, \cons \rangle$ underlying the 3-valued logical matrix ${\cal M}_e$, where the operations are defined as in Theorem \ref{SoundCompFI2}. Let $X$ be a non-empty set. Recall that a {\em triple over $X$}  is a map $r:X\fun{\bf 3}$. If $r$ is a triple over $X$ we write $r= \langle r_{ \oplus},r_{ \ominus},r_{ \odot}\rangle$
where  \  \
$r_{ \oplus}=r^{-1}(1)$, \ 
$r_{ \odot}=r^{-1}(\frac{1}{2})$, \
$r_{ \ominus}=r^{-1}(0)$.
As usual, we denote by ${\bf 3}^X$ the set of all triples over $X$. Clearly, the set ${\bf 3}^X$ of triples inherits the algebraic structure of $\mathsf{3}$, where the  operations are defined pointwise. In a sense, ${\bf 3}^X$ generalizes the power-set ${\bf 2}^X$ (seen as a Boolean algebra), and so ${\bf 3}^X$ is a kind of 3-valued power-set, endowed with the 3-valued algebraic structure over $\{\wedge, \vee, \to, \neg, \circ \}$ inherited from the 3-valued matrix of \ciore.

\begin{prop}(\cite{ConGomFil})\label{lfi1al}
Let $r = (r_\oplus,r_\ominus,r_\odot)$ and  $u = (u_\oplus,u_\ominus,u_\odot)$ be two triples over $X$. Then:\\[1mm]
(i) $r\wedge u = ( ( r_\oplus \cap u_\oplus)\cup( r_\oplus\cap  u_\odot)\cup( r_\odot\cap  u_\oplus),  r_\ominus\cup  u_\ominus,  r_\odot \cap  u_\odot)$,\\[1mm]
(ii) $ r \vee u  = (  r_\oplus\cup u_\oplus \cup (r_\odot\cap u_\ominus)\cup (r_\ominus \cap u_\odot), r_\ominus\cap u_\ominus, r_\odot\cap u_\odot)$,\\[1mm]
(iii) $ r\imp u = ( r_\ominus\cup u_\oplus \cup (r_\oplus\cap u_\odot), (r_\oplus\cup r_\odot)\cap u_\ominus, r_\odot\cap u_\odot) $,\\[1mm]
(iv) $\no r = (r_\ominus,r_\oplus,r_\odot)$,\\[1mm]
(v) $\cons r = (r_\oplus \cup r_\ominus, r_\odot,\emptyset)$.
\end{prop}

It is not difficult to see that, if $r \in {\bf 3}^X$, then: (1)~$r_\oplus\cup r_\ominus\cup r_\odot=X$; and (2)~$r_\oplus\cap r_\ominus=r_\oplus\cap r_\odot=r_\ominus\cap r_\odot=\emptyset$.
Conversely, if $r_1$, $r_2$ and $r_3$ are subsets of $X\neq\emptyset$  such that (1)~$r_1\cup r_2\cup r_3=X$ and (2)~$r_1\cap r_2=r_1 \cap r_3=r_2\cap r_3=\emptyset$, then there exists a unique $r\in {\bf 3}^X$  such that $r=( r_1, r_2, r_3)$.\\
A {\em partial (or pragmatic) structure for \qciore\  over the signature $ \Theta$}  is an ordered pair $\ff A=\langle A,(\cdot)^\ff A\rangle$ 
where $A\neq\emptyset$ and $(\cdot)^\ff A$ is a function such that

\begin{itemize}
\item[-] for all $R\in\mm P_n$, 
$R^{\ff A}=\langle R^{\ff A}_{ \oplus},R^{\ff A}_{ \ominus},R^{\ff A}_{ \odot}\rangle$
is a triple  over $A^n$,
\item[-] $(\cdot)^\ff A$ is defined as usual over $\mm F$ and $\mm C$.
\end{itemize}


\noindent Let $A$ be a non-empty set.  An {\em assignment into $A$} is any map $s:{\mm V}={\mm V}_f\cup{\mm V}_b\rightarrow  A$. We denote by $S(A)$ the set of all assignments into $A$, i.e. $S(A)=A^\mm V$. If $\ff A$ is a partial structure over $\Theta$ with domain $A$ then an {\em assignment into $\ff A$} is any assignment into $A$. The set of all assignments into $\ff A$ will be denoted by $S(\ff A)$, i.e. $S(\ff A)=S(A)=A^\mm V$.

\begin{defi}\label{def-valterm}
Let $s \in S(\ff A)$. The {\em value of the term $t$ in $\ff A$ under the assignment $s$}, denoted by $t^\ff A[s]$, is defined inductively as follows:

\begin{itemize}
\item[-] if $t$ is $v$, for $v\in\mm V$, then $t^\ff A[s]:=s(v)$,
\item[-] if $t$ is $c$, for $c\in\mm C$, then $t^\ff A[s]:=c^\ff A$,
\item[-] if $t$ is $f(t_1,\ldots,t_n)$, for $f\in\mm F_n$ and terms $t_i$, then $t^\ff A[s]:=f^\ff A(\ter{1},\ldots,\ter{n})$.
\end{itemize}
\end{defi}

\begin{defi}(\cite{ConGomFil})\label{def-quant} 
Let $A$ be a non-empty set. Given a set $Z$ let $\wp(Z)_+$ be the set  of  non-empty subsets of $Z$. For $a\in {\mm V}_f$ let $\widehat{\pt a}:\wp(S(A)) \to  \wp(S(A))$ and $\widehat{\ex a}:\wp(S(A)) \to  \wp(S(A))$ be defined as follows, for every $Y \subseteq S(A)$:

$${\widehat{\pt a}(Y)}:=\{s\in S(A) \ : \  s_a^m\in Y\mbox{ for all }m\in A\},$$
$${\widehat{\ex a}(Y)}:=\{s\in S(A) \ : \  s_a^m\in Y\mbox{ for some }m\in A\}.$$
The functions $\widetilde{\forall}:\wp({\bf 3})_+ \to {\bf 3}$ and  $\widetilde{\exists}:\wp({\bf 3})_+ \to {\bf 3}$ are defined as follows, for every $\emptyset\neq Y \subseteq {\bf 3}$:
$$
\widetilde{\forall}(Y):=
\left\{
\begin{array}{ll}
\1
& \ \mbox{if }\ 
\1 \in Y, \ \0 \notin Y
\\[1mm]
\um
& \ \mbox{if }\ 
Y=\{\um\}
\\[1mm]
\0
& \ \mbox{if }\ 
\0 \in Y
\end{array}
\right.
\hspace{1.5cm}
\widetilde{\exists}(Y):=
\left\{
\begin{array}{ll}
\1
& \ \mbox{if }\ 
Y\neq\{\um\}, \ Y\neq\{\0\}
\\[1mm]
\um
& \ \mbox{if }\ 
Y=\{\um\}
\\[1mm]
\0
& \ \mbox{if }\ 
Y=\{\0\}
\end{array}
\right.
$$
\end{defi}

\begin{defi}(\cite{ConGomFil}) \label{defLFI2structure} A {\em \qciore-structure} over $\Theta$ is a pair
$\langle\ff A,\tripc{\cdot}\rangle$ such that $\ff A$ is a partial structure for \qciore\ and  $\tripc{\cdot}:\rr L(\Theta)\imp {\bf 3}^{S(\ff A)}$ is a map defined recursively as follows, for every $s \in S(\ff A)$:\\[1mm]
1. If $P(t_1,\ldots,t_n)$ is atomic, $\tripc{P(t_1,\ldots,t_n)}(s)=
P^\ff A(t_1^{\ff A}[s],\ldots,t_n^{\ff A}[s])$,\\[1mm]
2. $\tripc{\varphi\wedge\psi}=\tripc{\varphi}\wedge \tripc{\psi}$,\\[1mm]
3. $\tripc{\varphi\vee\psi}=\tripc{\varphi}\vee \tripc{\psi}$,\\[1mm]
4. $\tripc{\varphi\imp\psi}=\tripc{\varphi} \imp \tripc{\psi}$,\\[1mm]
5. $\tripc{\no \varphi}=\no \tripc{\varphi}$,\\[1mm]
6. $\tripc{\cons\varphi}=\cons\tripc{\varphi}$,\\[1mm] 
7. $\tripc{\pt x \varphi(x)}(s)=\widetilde{\forall}(\{\tripc{\varphi(a)}(s_a^m) \ : \ m\in A\})$, where $\varphi(x)$ is $\left(\varphi(a)\, \d \frac{a}{x}\right)$,\\[1mm] 
8. $\tripc{\ex x \varphi(x)}(s)=\widetilde{\exists}(\{\tripc{\varphi(a)}(s_a^m) \ : \ m\in A\})$, where $\varphi(x)$ is $\left(\varphi(a)\, \d \frac{a}{x}\right)$.
\end{defi}

\begin{prop}(\cite{ConGomFil}) \label{charLFI2structure} Let 
$\langle\ff A,\|\cdot\|^\ff A\rangle$ be a {\em \quci-structure} over a signature $\Theta$. Then, by denoting $\tripc{\varphi}$ as $\langle \mas\varphi, \menos\varphi, \por\varphi \rangle$, the following holds:
 \begin{itemize}
 
 \item[] If $P(t_1,\ldots,t_n)$ is atomic then, for every $\# \in \{\oplus,\ominus,\odot\}$ and $s \in S(\ff A)$:
$$
s\in \|P(t_1, \ldots,t_n)\|^{\mathfrak A}_\#
 \ \mbox{ iif } \ 
(t_1^{\ff A}[s],\ldots,t_n^{\ff A}[s])\in P^{\mathfrak A}_\#,$$
 
\item[]
$\tripc{\varphi\wedge\psi} \ := \ 
\left\langle
(\mas\varphi\cap\mas\psi)\cup
(\mas\varphi\cap\por\psi)\cup
(\por\varphi\cap\mas\psi),\right.$

$\phantom{\tripc{\varphi\wedge\psi} \ := \ } \ \
\left.\menos\varphi\cup\menos\psi,
\por\varphi\cap\por\psi
\right\rangle
$,\\

\item[]
$\tripc{\varphi\vee\psi} \ := \ 
\left\langle
\mas\varphi\cup\mas\psi\cup(\por\varphi\cap\menos\psi)\cup(\menos\varphi\cap\por\psi),\right.$

$\phantom{\tripc{\varphi\vee\psi} \ := \ } \ 
\left.\menos\varphi\cap\menos\psi, \por\varphi\cap\por\psi
\right\rangle
$,\\

\item[]
$\tripc{\varphi\imp\psi} \ := \ \left\langle
\menos\varphi\cup \mas\psi\cup(\mas\varphi\cap\por\psi), (\mas\varphi \cup
\por\varphi)\cap\menos\psi, \por\varphi\cap\por\psi
\right\rangle
$,\\

\item[]
$\tripc{\no \varphi}\ := \ 
\langle
\menos\varphi,
\mas\varphi,
\por\varphi
\rangle$,\\

\item[]
$\tripc{\cons\varphi}\ :=\ 
\langle
\mas\varphi\cup\menos\varphi,
\por\varphi,
\emptyset
\rangle$,\\

\item[]
$\tripc{\pt x \varphi(x)}\ := 
\left\langle 
\widehat{\ex a}(\mas{\varphi(a)})-\widehat{\ex a}(\menos\varphi(a)),
\widehat{\ex a}(\menos{\varphi(a)}),
\widehat{\pt a}(\por{\varphi(a)})
\right\rangle$, where $\varphi(x)$ is $\left(\varphi(a)\, \d \frac{a}{x}\right)$,\\

\item[]
$\tripc{\ex x \varphi(x)}\ := 
\left\langle 
S(\ff A)-\big(\widehat{\pt a}(\menos{\varphi(a)})\cup\widehat{\pt a}(\por{\varphi(a)})\big),\ 
\widehat{\pt a}(\menos{\varphi(a)}),\ 
\widehat{\pt a}(\por{\varphi(a)})
\right\rangle$, where $\varphi(x)$ is $\left(\varphi(a)\, \d \frac{a}{x}\right)$.
 \end{itemize}
\end{prop}

\noindent Observe that a partial structure for  \qciore\ over $\Theta$ determines a unique \qciore-structure over $\Theta$. Hence, both structures will be identified from now on. 
\begin{defi}(\cite{ConGomFil})  \label{defsat}
Let $\langle\ff A,\tripc{\cdot}\rangle$ be a \qciore-structure, let $\varphi$ be a formula and let $s\in\ese A$.  
We say that $s$ {\em satisfies} $\varphi$ in $\ff A$, denoted by $\ff A\models\varphi[s]$, if $s\in\mas\varphi\cup\por\varphi$. Besides,  
$\varphi$ is said to be {\em valid} in $\ff A$ (or that $\ff A$ validates $\varphi$), denoted by $\ff A\models_\qciore\varphi$ (or simply $\ff A\models\varphi$) , if  $\mas{\varphi}\cup\por{\varphi}\ =\ S(\ff A)$. A \qciore-structure $\ff A$ is a model of a set $\Gamma$ of formulas if $\ff A\models_\qciore\varphi$ for each $\varphi\in\Gamma$. Besides, given a set $\Gamma\cup \{ \varphi \}$ of formulas, we say that {\em $\varphi$ is a \qciore-consequence of $\Gamma$}, denoted by $\Gamma\models_\qciore\varphi$ if, for every
\qciore-structure $\ff A$, we have that: $\ff A\models_\qciore\psi$ for every $\psi\in \Gamma$ implies that $\ff A\models_\qciore\varphi$.
\end{defi}

\begin{theo} (\cite{ConGomFil}) [Soundness and Completeness] \label{soundQmbC} 
Let $\Gamma \cup \{\varphi\}$ be a set of formulas over $\Gamma$. 
$$\Gamma\vdash_{\bf QCiore}\varphi \, \mbox{ iff } \, \Gamma\models_\qciore\varphi$$
\end{theo}

\

\section{Gentzen-style proof systems for \ciore}\label{s3}

In this section, we introduce a cut-free sequent calculus for \ciore. Here, we strongly rely on the general method, described in \cite{Avron02}, to obtain cut-free systems for logics with matrix semantics. However, here we do not describe the full process (as it is described in \cite{MF} and \cite{LCMF}) carried out to obtain such systems. Instead, we present a sequent calculus and provide proofs for the soundness, completeness and cut-elimination theorem.

In \cite{Avron01}, it was developed a generic $n$-sequent system for any logic based on an $n$-valued (non-deterministic) matrix. The following is a 3-sequent system for \ciore\ which enjoys the cut-elimination property.

\begin{defi} The 3-sequent calculus $\bf{3-{\cal S}}Ciore$ is the system defined as follows: for  $\alpha, \beta\in Fm$, $\Omega$ and $\Omega'$ arbitrary sets of signed formulas

\begin{itemize}
\item {\bf Axioms:} $\{\0:\alpha, \um:\alpha, \1:\alpha\}$.
\item {\bf Structural rules:} Weakening: $ \displaystyle\frac{\Omega}{\Omega '}$ in case $\Omega\subseteq\Omega '$.
\item {\bf Logical rules:} For $i, j \in {\cal T}$
\end{itemize}

$$ \mbox{\rm($\vee_{ij})$ \, } \displaystyle \frac{\Omega, i:\alpha \hspace{0.5cm} \Omega, j:\beta} {\Omega, 1:\alpha\vee\beta} \,\mbox{ for } i=j=\1 \mbox{ or } i \neq j  \hspace{1.5cm}  \mbox{\rm($\vee_{\0\0})$ \, } \displaystyle \frac{\Omega, \0:\alpha  \hspace{0.5cm} \Omega, \0:\beta} {\Omega, \0:\alpha\vee\beta}$$
$$\mbox{\rm($\vee_{\um\um})$ \, } \displaystyle \frac{\Omega, \um :\alpha \hspace{0.5cm}\Omega, \um :\beta} {\Omega, \um :\alpha\vee\beta}$$

$$\mbox{\rm($\wedge_{ij})$ \, } \displaystyle \frac{\Omega, i:\alpha  \hspace{0.5cm} \Omega, j:\beta} {\Omega, \1:\alpha\wedge\beta} \mbox{ for } i=\1, j\neq \0 \, \mbox{ or } \, j=\1, i \neq \0$$

$$ \mbox{\rm($\wedge_{ij})$ \, } \displaystyle \frac{\Omega, i:\alpha\hspace{0.5cm} \Omega, j:\beta} {\Omega, \0:\alpha\wedge\beta} \mbox{ for } i=\0\, \mbox{ or } \, j=\0  \hspace{1.5cm} \mbox{\rm($\wedge_{ij})$ \, } \displaystyle \frac{\Omega, i:\alpha \hspace{0.5cm} \Omega, j:\beta} {\Omega, \um:\alpha\wedge\beta} \mbox{ for } \, i=j=\um$$

 $$   \mbox{\rm($\rightarrow_{ij})$ \, } \displaystyle \frac{\Omega, i:\alpha \hspace{0.5cm} \Omega, j:\beta} {\Omega, \1:\alpha\rightarrow\beta} \mbox{ for }\, i=\0 \, \mbox{ or } \, j=\1 \, \mbox{ or } \, (i=\1, j= \um) $$

$$ \mbox{\rm($\rightarrow_{i\0})$ \, } \displaystyle \frac{\Omega, i:\alpha \hspace{0.5cm} \Omega, \0:\beta} {\Omega, \0:\alpha\rightarrow\beta} \mbox{ for }\, i \neq \0  \hspace{1.5cm} \mbox{\rm($\rightarrow_{\um\um})$ \, } \displaystyle \frac{\Omega, \um:\alpha \hspace{0.5cm} \Omega, \um:\beta} {\Omega, \um:\alpha\rightarrow\beta}\mbox{ for } \,  i=j=\um$$

$$ \mbox{\rm($\neg_{\0})$ \, } \displaystyle \frac{\Omega, \0:\alpha} {\Omega, \1:\neg \alpha} \hspace{0.5cm}
 \mbox{\rm($\neg_{\um})$ \, } \displaystyle \frac{\Omega, \um:\alpha} {\Omega, \um:\neg \alpha} \hspace{0.5cm}
 \mbox{\rm($\neg_{\1})$ \, } \displaystyle \frac{\Omega, \1:\alpha} {\Omega, \0:\neg \alpha}$$
$$ \mbox{\rm($\circ_{\0})$ \, } \displaystyle \frac{\Omega, \0:\alpha} {\Omega, \1:\circ \alpha}\hspace{0.5cm}
\mbox{\rm($\circ_{\um})$ \, } \displaystyle \frac{\Omega, \um:\alpha} {\Omega, \0:\circ \alpha} \hspace{0.5cm}
 \mbox{\rm($\circ_{\1})$ \, } \displaystyle \frac{\Omega, \1:\alpha} {\Omega, \1:\circ \alpha}$$
\end{defi}

\

\begin{prop}
\begin{itemize}
\item[]
\item[\rm (i)] $\bf{3-{\cal S}Ciore}$ is sound and complete w.r.t. the matrix ${\cal M}_{e}$,
\item[\rm (ii)] The cut rule is admissible in  $\bf{3-{\cal S}Ciore}$.
\end{itemize}
\end{prop}
\begin{dem} By construction (see \cite{Avron01}).
\end{dem}

\

The above system can be translated to an ordinary two-sided sequent system providing that the language of \ciore\ is {\em sufficiently expressive} (cf. \cite{Avron02}). Recall that a language $\mathscr{L}$ is sufficiently expressive for the matrix ${\cal M}$, with set of truth values ${\cal T}=\{t_1,\dots,t_n\}$, iff for any $i$, $0\leq i \leq n-1$ there exist natural numbers $l_i,m_i$  and formulas $\alpha_{j}^{i}, \beta_{k}^{i}$ that have $p$ as their only propositional variable, for $1\leq j\leq l_i$ and $1\leq k\leq m_i$ such that for any valuation $v$, the following conditions hold:\\[1.5mm]
(i) $\alpha_{1}^{i}=p$ if $t_i\in {\cal N}$ \, and \, $\beta_{1}^{i}=p$ if $t_i\in {\cal D}$, \\[1.5mm]
(ii) For any $\varphi \in Fm$ and $t_i\in {\cal T}$
\begin{center}
\begin{tabular}{ccl}
$v(\varphi)=t_i$ &\, $\Leftrightarrow$ \, & $v(\alpha_{1}^{i}[p/\varphi]), \dots, v(\alpha_{l_i}^{i}[p/\varphi]) \in {\cal N}$ \, and \,  \\
 & & $v(\beta_{1}^{i}[p/\varphi]), \dots, v(\alpha_{m_i}^{i}[p/\varphi]) \in {\cal D}$ 
\end{tabular}
\end{center}
\noindent where $\alpha_{j}^{i}[p/\varphi]$ ($\beta_{k}^{i}[p/\varphi]$) is the formula obtained by the substitution of $p$ by $\varphi$ in $\alpha_{j}^{i}$ ($\beta_{k}^{i}$); and $\mm N = \mm T\setminus \mm D$.

\

The language of \ciore\ is sufficiently expressive for the semantics determined by the matrix ${\cal M}_e$. Indeed, if $v:Fm \to M_e$ is a valuation and $\alpha\in Fm$ is an arbitrary formula, then we have  that
$$v(\alpha)=0 \Longleftrightarrow v(\alpha) \in {\cal N}. $$
$$v(\alpha)=\frac{1}{2} \Longleftrightarrow v(\alpha) \in {\cal D} \, \mbox{ and } \, v(\neg \alpha) \in {\cal D}.$$
$$v(\alpha)=1 \Longleftrightarrow v(\neg\alpha) \in {\cal N}\in {\cal D} \, \mbox{ and } \, v( \alpha) \in {\cal D}. $$

After the translation, the system obtained needs (in general) to be streamlined to reduce it to a more
compact form. The details of the translation process of $\bf{3-{\cal S}Ciore}$ to an ordinary two-sided sequent system is spare here. However we shall mention the three general streamlining principles from \cite{Avron01} since we shall use them in what follows.  \\
Of these three, the first and the third decrease the number of rules (which is our main measure of complexity), while the second simplifies a rule by decreasing the number of its premises (since the third rule increases this number, its application is often followed by applications of the first two). \\[2mm]
Recall that a rule (r) is context-free if whenever \, $\displaystyle\frac{\Gamma_1\Rightarrow\Delta_1 \, \dots \, \Gamma_k\Rightarrow\Delta_k}{\Gamma\Rightarrow\Delta}$ \,  is a valid application of (r), and $\Gamma'$ and $\Delta'$ are sets of formulas, then 
$$\displaystyle\frac{\Gamma_1,\Gamma'\Rightarrow\Delta_1, \Delta' \, \dots \, \Gamma_k, \Gamma' \Rightarrow\Delta_k, \Delta'}{\Gamma, \Gamma' \Rightarrow\Delta, \Delta'}$$
 is also a valid application of (r). The three streamlining principles above-mentioned are: 

\

\noindent {\bf Principle 1.} If a rule in is derivable from other rules, it can be deleted.

\

\noindent {\bf Principle 2.} If $\frac{S}{\Sigma}$ (where S is a set of premises) is a rule, $S'$ is a subset of $S$ and $\frac{S'}{\Sigma}$ is derivable (perhaps using cuts), then $\frac{S}{\Sigma}$ can be replaced by $\frac{S'}{\Sigma}$. In particular: if $\frac{S}{\Sigma}$  is a rule, $\Gamma\Rightarrow \Delta \in S$, and $\Gamma\Rightarrow \Delta$ is derivable from $S\setminus \{\Gamma\Rightarrow \Delta\}$, then $\frac{S}{\Sigma}$ can be replaced with $\frac{S\setminus \{\Gamma\Rightarrow \Delta\}}{\Sigma}$. Two very simple, but
quite useful cases of this are when $\Gamma\Rightarrow \Delta$ is subsumed by an axiom or by some sequent in $\Gamma\Rightarrow \Delta$.

\

\noindent {\bf Principle 3.} If in a given sequent system $G$ we have the rules $\displaystyle \frac{\{\Gamma_{i}\Rightarrow\Delta_{i}\}_{\small 1 \leq i \leq n}} {\Gamma\Rightarrow\Delta}$ , $\displaystyle \frac{\{\Pi_{j}\Rightarrow\Theta_{j}\}_{1 \leq j\leq k} }{\Gamma\Rightarrow\Delta}$ and both are context--free, then we can replace these two rules by the new rule   
$$\displaystyle \frac{\{\Gamma_{i},\Pi_{j}\Rightarrow\Delta_{i},\Theta_{j}\}_{\small  1\leq i \leq n,  1 \leq j \leq k }} {\Gamma\Rightarrow\Delta}$$

\begin{rem}\label{rem1} It is worth mentioning that if the cut-rule is admissible in a given ordinary sequent system then the same is true for the systems obtained from it by using the above streamlining principles, even if cuts are used in applications of  Principle 2  \cite[pg. 47]{Avron02}.
\end{rem}

\noindent The following is the sequent system obtained and which we call \gciore .

\begin{center}
{\bf Axioms}
$$\alpha \Rightarrow \alpha$$
{\bf Structural rules}
$$\mbox{\rm(w$\Rightarrow$) \, }  \displaystyle \frac{\Gamma\Rightarrow\Delta} {\Gamma,\alpha\Rightarrow\Delta} \hspace{2cm}  \mbox{\rm($\Rightarrow$w) }  \displaystyle \frac{\Gamma\Rightarrow\Delta} {\Gamma\Rightarrow\Delta, \alpha}$$
{\bf Logical rules}
$$ \mbox{\rm($\vee\Rightarrow)$ \, } 
   \displaystyle \frac{\Gamma,\alpha\Rightarrow\Delta \hspace{.5cm} \Gamma,\beta\Rightarrow\Delta} {\Gamma,\alpha\vee\beta\Rightarrow\Delta} \hspace{1.5cm} \mbox{\rm($\Rightarrow\vee)$ \, } \displaystyle \frac{\Gamma\Rightarrow\Delta,\alpha, \beta} {\Gamma\Rightarrow\Delta,\alpha\vee\beta}$$
$$ \mbox{\rm($\neg\vee\Rightarrow)$ \, }   \displaystyle \frac{\Gamma,\alpha,\neg\alpha, \beta, \neg\beta \Rightarrow\Delta \hspace{0.5cm} \Gamma,\neg\alpha, \neg\beta \Rightarrow\Delta, \alpha, \beta} {\Gamma, \neg(\alpha\vee\beta)\Rightarrow\Delta}$$
$$ \mbox{\rm($\Rightarrow\neg\vee)$ \, } 
   \displaystyle \frac{\Gamma\Rightarrow\Delta,\alpha \hspace{.5cm}\Gamma\Rightarrow\Delta,\neg\alpha \hspace{.5cm} \Gamma\Rightarrow\Delta,\beta\hspace{.5cm} \Gamma\Rightarrow\Delta,\neg\beta} {\Gamma\Rightarrow\Delta,\neg(\alpha\vee\beta)}$$

$$ \mbox{\rm($\wedge\Rightarrow)$ \, } 
   \displaystyle \frac{\Gamma,\alpha,\beta\Rightarrow\Delta} {\Gamma,\alpha\wedge\beta\Rightarrow\Delta} \hspace{1.5cm} \mbox{\rm($\Rightarrow\wedge)$ \, } \displaystyle \frac{\Gamma\Rightarrow\Delta,\alpha \hspace{.5cm}\Gamma\Rightarrow\Delta,\beta} {\Gamma\Rightarrow\Delta,\alpha\wedge\beta}$$
$$ \mbox{\rm($\neg\wedge\Rightarrow)$ \, } 
  \displaystyle \frac{\Gamma\Rightarrow\Delta, \alpha, \beta \hspace{0.5cm} \Gamma,\neg\alpha \Rightarrow\Delta, \alpha \hspace{0.5cm} \Gamma,\neg\beta \Rightarrow\Delta, \beta \hspace{0.5cm} \Gamma,\neg\alpha, \neg\beta \Rightarrow\Delta} {\Gamma, \neg(\alpha\wedge\beta)\Rightarrow\Delta} $$
 $$ \mbox{\rm($\Rightarrow\neg\wedge)$ \, } 
   \displaystyle \frac{\Gamma\Rightarrow\Delta,\alpha \hspace{.5cm}\Gamma\Rightarrow\Delta,\neg\alpha \hspace{.5cm} \Gamma\Rightarrow\Delta,\beta\hspace{.5cm} \Gamma\Rightarrow\Delta,\neg\beta} {\Gamma\Rightarrow\Delta,\neg(\alpha\wedge\beta)}$$

$$ \mbox{\rm($\rightarrow\Rightarrow)$ \, } 
   \displaystyle \frac{\Gamma\Rightarrow\Delta,\alpha \hspace{.5cm} \Gamma,\beta\Rightarrow\Delta} {\Gamma,\alpha\rightarrow\beta\Rightarrow\Delta} \hspace{1cm} \mbox{\rm($\Rightarrow\rightarrow)$ \, } \displaystyle \frac{\Gamma,\alpha \Rightarrow\Delta,\beta} {\Gamma\Rightarrow\Delta,\alpha\rightarrow\beta} $$
$$ \mbox{\rm($\neg\rightarrow\Rightarrow)$ \, }     \displaystyle \frac{\Gamma,\alpha,\neg\beta\Rightarrow\Delta, \beta \hspace{0.5cm} \Gamma,\alpha, \neg\alpha, \neg\beta \Rightarrow\Delta} {\Gamma, \neg(\alpha\rightarrow\beta)\Rightarrow\Delta} $$
 $$ \mbox{\rm($\Rightarrow\neg\rightarrow)$ \, } 
   \displaystyle \frac{\Gamma\Rightarrow\Delta,\alpha \hspace{.5cm}\Gamma\Rightarrow\Delta,\neg\alpha \hspace{.5cm} \Gamma\Rightarrow\Delta,\beta\hspace{.5cm} \Gamma\Rightarrow\Delta,\neg\beta} {\Gamma\Rightarrow\Delta,\neg(\alpha\rightarrow\beta)}$$

$$ \mbox{\rm($\Rightarrow\neg)$ \, } \displaystyle \frac{\Gamma,\alpha\Rightarrow\Delta} {\Gamma\Rightarrow\Delta,\neg\alpha}\hspace{1.3cm} \mbox{\rm($\neg\neg\Rightarrow)$ \, } \displaystyle \frac{\Gamma,\alpha\Rightarrow\Delta} {\Gamma,\neg\neg\alpha\Rightarrow\Delta} \hspace{1.3cm} \mbox{\rm($\Rightarrow\neg\neg)$ \, }   \displaystyle \frac{\Gamma\Rightarrow\Delta,\alpha} {\Gamma\Rightarrow\Delta,\neg\neg\alpha}$$

$$ \mbox{\rm($\circ\Rightarrow)$ \, } \displaystyle \frac{\Gamma\Rightarrow\Delta,\alpha \hspace{.5cm} \Gamma\Rightarrow\Delta,\neg\alpha} {\Gamma,\circ\alpha\Rightarrow\Delta} \hspace{1cm} \mbox{\rm($\Rightarrow\circ)$ \, } \displaystyle \frac{\Gamma,\alpha, \neg\alpha \Rightarrow\Delta} {\Gamma\Rightarrow\Delta, \circ\alpha} $$
$$  \mbox{\rm($\neg\circ\Rightarrow)$ \, } \displaystyle \frac{\Gamma, \alpha, \neg\alpha \Rightarrow\Delta} {\Gamma, \neg\circ\alpha\Rightarrow\Delta}$$
\end{center}

\

\noindent As usual, we write $\vdash_{\qciore}\Gamma\Rightarrow\Delta$ to indicate that there is a proof of $\Gamma\Rightarrow\Delta$ in \gciore . Next, we show some cut-free proofs of important theorems of \ciore\ which were stated in \cite{ConGomFil}. As usual, we write $\alpha\leftrightarrow\beta$ as an abbreviation of $(\alpha\to\beta)\wedge(\beta\to\alpha)$.

\begin{theo}
The following sequents are (cut-free) provable in \qciore.
\begin{enumerate}[{\rm (i)}]
\item $\Rightarrow(\alpha\wedge\neg\alpha)\leftrightarrow\neg\circ\alpha$,
\item $\Rightarrow\circ\circ\alpha$,
\item $\Rightarrow\circ\alpha\leftrightarrow\circ\neg\alpha$,
\item $\Rightarrow((\alpha\wedge\neg\alpha)\wedge(\beta\wedge\neg\beta))\leftrightarrow((\alpha\wedge\beta)\wedge\neg(\alpha\wedge\beta))$,
\item $\Rightarrow((\alpha\wedge\neg\alpha)\wedge(\beta\wedge\neg\beta))\leftrightarrow((\alpha\vee\beta)\wedge\neg(\alpha\vee\beta))$,
\item $\Rightarrow((\alpha\wedge\neg\alpha)\wedge(\beta\wedge\neg\beta))\leftrightarrow((\alpha\rightarrow\beta)\wedge\neg(\alpha\rightarrow\beta))$,
\item $\Rightarrow (\circ\alpha\vee\circ\beta)\sii\circ(\alpha\wedge\beta)$,
\item $\Rightarrow (\circ\alpha\vee\circ\beta)\sii\circ(\alpha\vee\beta)$,
\item $\Rightarrow (\circ\alpha\vee\circ\beta)\sii\circ(\alpha\imp\beta)$.
\end{enumerate}
\end{theo}

\begin{dem} We shall only prove items (i), (ii), (iii) and (vi). The rest are analogous.\\
(i) It is consequence of the following proofs.

\

\AxiomC{$\alpha \Rightarrow \alpha$}
\LeftLabel{\small(w$\Rightarrow$)}
\UnaryInfC{$\neg \alpha, \alpha\Rightarrow\alpha$}
\AxiomC{$\neg\alpha \Rightarrow \neg\alpha$}
\LeftLabel{\small(w$\Rightarrow$)}
\UnaryInfC{$\neg \alpha, \alpha\Rightarrow\neg\alpha$}
\LeftLabel{\small( $\circ\Rightarrow$ )}
\BinaryInfC{$\neg \alpha, \alpha, \circ\alpha\Rightarrow $}
\LeftLabel{\small( $\Rightarrow\neg$ )}
\UnaryInfC{$\neg \alpha, \alpha\Rightarrow \neg\circ\alpha $}
\LeftLabel{\small( $\wedge\Rightarrow$ )}
\UnaryInfC{$\neg \alpha\wedge \alpha\Rightarrow \neg\circ\alpha$}
\LeftLabel{\small( $\Rightarrow\rightarrow$ )}
\UnaryInfC{$\Rightarrow(\neg \alpha\wedge \alpha)\rightarrow \neg\circ\alpha$}
\DisplayProof \hspace{.8cm}
\AxiomC{$\alpha \Rightarrow \alpha$}
\LeftLabel{\small(w$\Rightarrow$)}
\UnaryInfC{$\neg \alpha, \alpha\Rightarrow\alpha$}
\LeftLabel{\small($\neg\circ\Rightarrow)$}
\UnaryInfC{$\neg \circ\alpha\Rightarrow\alpha$}
\AxiomC{$\neg\alpha \Rightarrow \neg\alpha$}
\LeftLabel{\small(w$\Rightarrow$)}
\UnaryInfC{$\neg \alpha, \alpha\Rightarrow\neg\alpha$}
\LeftLabel{\small($\neg\circ\Rightarrow)$}
\UnaryInfC{$\neg \circ\alpha\Rightarrow\neg\alpha$}
\LeftLabel{\small( $\Rightarrow\wedge$ )}
\BinaryInfC{$\neg \circ\alpha\Rightarrow (\alpha\wedge\neg\alpha) $}
\LeftLabel{\small( $\Rightarrow\rightarrow$ )}
\UnaryInfC{$\Rightarrow\neg \circ\alpha\rightarrow (\alpha\wedge\neg\alpha)  $}
\DisplayProof

\

\noindent  (ii)
 
\begin{prooftree}
\AxiomC{$\alpha\Rightarrow\alpha$}
\LeftLabel{\small(w$\Rightarrow$)}
\UnaryInfC{$\neg\alpha, \alpha\Rightarrow\alpha$}
\AxiomC{$\neg\alpha\Rightarrow\neg\alpha$}
\LeftLabel{\small(w$\Rightarrow$)}
\UnaryInfC{$\neg\alpha, \alpha\Rightarrow\neg\alpha$}
\LeftLabel{\small( $\circ\Rightarrow$ )}
\BinaryInfC{$\circ\alpha, \alpha, \neg\alpha\Rightarrow$}
\LeftLabel{\small( $\neg\circ\Rightarrow)$}
\UnaryInfC{$\circ\alpha, \neg\circ\alpha\Rightarrow $}
\LeftLabel{\small($\Rightarrow\circ)$}
\UnaryInfC{$\Rightarrow\circ\circ\alpha $}
\end{prooftree}

\

\noindent (iii)

\

\AxiomC{$\alpha\Rightarrow\alpha$}
\LeftLabel{\small($\neg\neg\Rightarrow$)}
\UnaryInfC{$\neg\neg\alpha\Rightarrow\alpha$}
\LeftLabel{\small(w$\Rightarrow$)}
\UnaryInfC{$\neg\neg\alpha, \neg\alpha\Rightarrow\alpha$}
\LeftLabel{\small($\Rightarrow\circ)$}
\UnaryInfC{$\Rightarrow\alpha, \circ\neg\alpha$}
\AxiomC{$\neg\alpha\Rightarrow\neg\alpha$}
\LeftLabel{\small(w$\Rightarrow$)}
\UnaryInfC{$\neg\neg\alpha, \neg\alpha\Rightarrow\neg\alpha$}
\LeftLabel{\small($\Rightarrow\circ)$}
\UnaryInfC{$\Rightarrow\neg\alpha, \circ\neg\alpha$}
\LeftLabel{\small( $\circ\Rightarrow$ )}
\BinaryInfC{$\circ\alpha\Rightarrow\circ\neg\alpha$}
\LeftLabel{\small($\Rightarrow\rightarrow)$}
\UnaryInfC{$\Rightarrow\circ\alpha\rightarrow\circ\neg\alpha$}
\DisplayProof \hspace{.8cm}
\AxiomC{$\alpha\Rightarrow\alpha$}
\AxiomC{$\alpha\Rightarrow\alpha$}
\LeftLabel{\small($\Rightarrow\neg\neg$)}
\UnaryInfC{$\alpha\Rightarrow\neg\neg\alpha$}
\LeftLabel{\small( $\circ\Rightarrow$ )}
\BinaryInfC{$\circ\neg\alpha, \alpha\Rightarrow$}
\LeftLabel{\small(w$\Rightarrow$)}
\UnaryInfC{$\circ\neg\alpha, \alpha, \neg\alpha\Rightarrow$}
\LeftLabel{\small($\Rightarrow\circ)$}
\UnaryInfC{$\circ\neg\alpha\Rightarrow\circ\alpha$}
\LeftLabel{\small($\Rightarrow\rightarrow)$}
\UnaryInfC{$\Rightarrow\circ\neg\alpha\rightarrow\circ\alpha$}
\DisplayProof

\

\

\noindent (vi) \,   One can easily find a derivation for $((\alpha\wedge\neg\alpha)\wedge(\beta\wedge\neg\beta))\rightarrow(\alpha\rightarrow\beta)$. On the other hand, a derivation for 
$((\alpha\wedge\neg\alpha)\wedge(\beta\wedge\neg\beta))\rightarrow\neg(\alpha\rightarrow\beta))$ is the following:

\begin{prooftree}\small
\AxiomC{$\alpha \Rightarrow \alpha$}
\UnaryInfC{$\neg \alpha, \alpha, \beta, \neg\beta \Rightarrow\alpha$}
\AxiomC{$\neg\alpha \Rightarrow \neg\alpha$}
\UnaryInfC{$\neg \alpha, \alpha, \beta, \neg\beta \Rightarrow\neg\alpha$}
\AxiomC{$\beta \Rightarrow \beta$}
\UnaryInfC{$\neg \alpha, \alpha, \beta, \neg\beta \Rightarrow\alpha$}
\AxiomC{$\neg\beta \Rightarrow \neg\beta$}
\UnaryInfC{$\neg \alpha, \alpha, \beta, \neg\beta \Rightarrow\neg\beta$}
\LeftLabel{\small( $\Rightarrow\neg\to$ )}
\QuaternaryInfC{$\neg \alpha, \alpha, \beta, \neg\beta \Rightarrow\neg(\alpha\to\beta) $}
\LeftLabel{\small( $\wedge\Rightarrow$'s )}
\UnaryInfC{$\alpha\wedge \neg\alpha, \beta\wedge\neg\beta \Rightarrow\neg(\alpha\to\beta) $}
\LeftLabel{\small( $\wedge\Rightarrow$ )}
\UnaryInfC{$(\alpha\wedge \neg\alpha)\wedge(\beta\wedge\neg\beta)\Rightarrow\neg(\alpha\to\beta) $}
\end{prooftree}

\noindent Then, by using ($\Rightarrow\wedge$), we have a derivation for  
$$\Rightarrow((\alpha\wedge\neg\alpha)\wedge(\beta\wedge\neg\beta))\rightarrow((\alpha\rightarrow\beta)\wedge\neg(\alpha\rightarrow\beta)).$$
In a similar way we show that 
$$\Rightarrow((\alpha\rightarrow\beta)\wedge\neg(\alpha\rightarrow\beta))\to((\alpha\wedge\neg\alpha)\wedge(\beta\wedge\neg\beta))$$
is derivable. 
\end{dem}

\begin{rem} Items (iv), (v) and (vi) say that, similarly to what happens with consistency, {\em contradictoriness} also propagates in \ciore .
\end{rem}

\section{Soundness, completeness and cut-elimination}\label{s4}

In this section, we show that the deductive system found in the previous section is sound and complete with respect to the matrix ${\cal M}_{e}$. Furthermore, we show that this sequent calculus enjoys the cut-elimination property.

\begin{theo}[Soundness] If \, $\vdash_{\gciore}\Gamma\Rightarrow\Delta$ \, then \, $\models_{{\cal M}_{e}}\Gamma\Rightarrow\Delta$.

\end{theo}
\begin{dem} It is enough to check that the axiom is valid and that every (structural and logic) rule preserves validity. The axiom is valid trivially, by definition of validity. Next, we shall see that ($\Rightarrow\neg\wedge$) is valid, the proof for the rest of the rules is similar. 
 $$ \mbox{\rm($\Rightarrow\neg\wedge)$ \, } 
   \displaystyle \frac{\Gamma\Rightarrow\Delta,\alpha \hspace{0.3cm}\Gamma\Rightarrow\Delta,\neg\alpha \hspace{0.3cm} \Gamma\Rightarrow\Delta,\beta \hspace{0.3cm} \Gamma\Rightarrow\Delta,\neg\beta} {\Gamma\Rightarrow\Delta,\neg(\alpha\wedge\beta)}$$
Suppose that $\models_{{\cal M}_{e}}\Gamma\Rightarrow\Delta,\alpha$, $\models_{{\cal M}_{e}}\Gamma\Rightarrow\Delta,\neg\alpha$, $\models_{{\cal M}_{e}}\Gamma\Rightarrow\Delta,\beta$ and $\models_{{\cal M}_{e}}\Gamma\Rightarrow\Delta,\neg\beta$. Let us see that  $\models_{{\cal M}_{e}}\Gamma\Rightarrow\Delta,\neg(\alpha\wedge\beta)$.
Let $v$ be a valuation. Since $v$ satisfies $\Gamma\Rightarrow\Delta,\alpha$, $\Gamma\Rightarrow\Delta,\neg\alpha$, $\Gamma\Rightarrow\Delta,\beta$ and $\Gamma\Rightarrow\Delta,\neg\beta$, we have to analize the following cases:\\
\underline{Case 1}: There exists $\gamma\in\Gamma$ such that $v(\gamma)\in{\cal N}$ or $\delta\in\Delta$ such that $v(\delta)\in{\cal D}$. In this case it is clear that $v$ satisfies $\Gamma\Rightarrow\Delta,\neg(\alpha\wedge\beta)$.\\
\underline{Case 2}: $v(\alpha)\in{\cal D}$, $v(\neg\alpha)\in{\cal D}$, $v(\beta)\in{\cal D}$ and $v(\neg\beta)\in{\cal D}$ then necessarily $v(\alpha)=\um$ and $v(\beta)=\um$. Therefore, $v(\alpha\wedge\beta)=\um$, by the table of $\hat\wedge$, and then $v(\neg(\alpha\wedge\beta))\in{\cal D}$. So, $v$ satisfies $\Gamma\Rightarrow\Delta,\neg(\alpha\wedge\beta)$.
\end{dem}

\

The following is an auxiliary result.

\begin{prop}\label{propaux} The following rules are derivable in \gciore .\\[3mm]
\begin{tabular}{ll}
{\rm (i)} \, $\mbox{\rm($\Rightarrow\neg\vee)'$ \, } \displaystyle \frac{\Gamma,\alpha \Rightarrow\Delta \hspace{.5cm}\Gamma, \beta\Rightarrow\Delta} {\Gamma\Rightarrow\Delta,\neg(\alpha\vee\beta)}$ & {\rm (ii)} \, $\mbox{\rm($\Rightarrow\neg\wedge)'$ \, } \displaystyle \frac{\Gamma,\alpha, \beta \Rightarrow\Delta} {\Gamma\Rightarrow\Delta,\neg(\alpha\wedge\beta)}$ \\[3mm]
{\rm (iii)} \, $\mbox{\rm($\Rightarrow\neg\rightarrow)'$ \, } \displaystyle \frac{\Gamma\Rightarrow\Delta, \alpha  \hspace{.5cm}\Gamma, \beta\Rightarrow\Delta} {\Gamma\Rightarrow\Delta,\neg(\alpha\rightarrow\beta)}$ & {\rm (iv)} \,  $\mbox{\rm($\Rightarrow\neg)'$ \, } \displaystyle \frac{\Gamma\Rightarrow\Delta,\alpha \,\,\,\,\,\,\Gamma\Rightarrow\Delta, \neg\alpha} {\Gamma\Rightarrow\Delta, \neg\alpha}$ 
\end{tabular}
\end{prop}

\begin{dem}

\begin{tabular}{clcl}
(i) &
\AxiomC{$\Gamma, \alpha \Rightarrow\Delta $}
\LeftLabel{\small($\Rightarrow \neg $)}
\AxiomC{$\Gamma, \beta \Rightarrow \Delta$}
\LeftLabel{\small($\neg\neg \Rightarrow$)}
\LeftLabel{\small($\vee\Rightarrow$)}
\BinaryInfC{$\Gamma,\alpha\vee\beta\Rightarrow\Delta$}
\LeftLabel{\small($\Rightarrow\neg$)}
\UnaryInfC{$\Gamma\Rightarrow\Delta, \neg(\alpha\vee\beta)$}
\DisplayProof \hspace{3cm}& (ii) & 
\AxiomC{$\Gamma, \alpha, \beta \Rightarrow\Delta $}
\LeftLabel{\small($\wedge\Rightarrow$)}
\UnaryInfC{$\Gamma, \alpha\wedge\beta\Rightarrow\Delta$}
\LeftLabel{\small($\Rightarrow\neg$)}
\UnaryInfC{$\Gamma\Rightarrow\neg(\alpha\wedge\beta)$}
\DisplayProof \\
& & & \\
(iii) & 
\AxiomC{$\Gamma\Rightarrow\Delta, \alpha $}
\AxiomC{$\Gamma, \beta \Rightarrow \Delta$}
\LeftLabel{\small($\rightarrow\Rightarrow$)}
\BinaryInfC{$\Gamma,\alpha\rightarrow\beta\Rightarrow\Delta$}
\LeftLabel{\small($\Rightarrow\neg$)}
\UnaryInfC{$\Gamma\Rightarrow\Delta \neg(\alpha\rightarrow\beta)$}
\DisplayProof & (iv) & Immediate.  \\
 & & &
\end{tabular}

\end{dem}

\

Let us consider the sequent calculus obtained from \gciore\ by adding the rules {\rm($\Rightarrow\neg\vee)'$}, {\rm($\Rightarrow\neg\wedge)'$}, {\rm($\Rightarrow\neg\rightarrow)'$} and {\rm($\Rightarrow\neg)'$} (as it is usual, we denote such system by $ \gciore\ + \mbox{\rm($\Rightarrow\neg\vee)'$} + \mbox{\rm($\Rightarrow\neg\wedge)'$} + \mbox{\rm($\Rightarrow\neg\rightarrow)'$} + \mbox{\rm($\Rightarrow\neg)'$}$). 
It is clear (see Proposition \ref{propaux}) that a sequent $\Gamma\Rightarrow\Delta$ is provable in \gciore\ iff  it is provable in  $ \gciore\ + \mbox{\rm($\Rightarrow\neg\vee)'$} + \mbox{\rm($\Rightarrow\neg\wedge)'$} + \mbox{\rm($\Rightarrow\neg\rightarrow)'$} + \mbox{\rm($\Rightarrow\neg)'$}$. Then, we modify this new system using only the principle 3, in the following  way: 
from $(\Rightarrow\neg\vee)$, $(\Rightarrow\neg\vee)'$ and Principle 3,  we obtain 
$$\mbox{\rm \small ($\Rightarrow\neg\vee)''$ \, } \small \frac{\Gamma,\alpha\Rightarrow\Delta, \neg\alpha  \hspace{.3cm}\Gamma,\alpha \Rightarrow\Delta, \beta \hspace{.3cm} \Gamma,\alpha \Rightarrow\Delta, \neg\beta \hspace{.3cm} \Gamma,\beta \Rightarrow\Delta, \alpha \hspace{.3cm} \Gamma,\beta \Rightarrow\Delta, \neg\alpha \hspace{.3cm} \Gamma,\beta \Rightarrow\Delta, \neg\beta} {\Gamma\Rightarrow\Delta, \neg(\alpha\vee\beta)}$$
From $(\Rightarrow\neg\wedge)$, $(\Rightarrow\neg\wedge)'$ and Principle 3, we obtain 
$$\mbox{\rm($\Rightarrow\neg\wedge)''$ \, } \displaystyle \frac{\Gamma,\alpha, \beta\Rightarrow\Delta, \neg\alpha  \hspace{.3cm}\Gamma,\alpha, \beta \Rightarrow\Delta, \neg\beta} {\Gamma\Rightarrow\Delta, \neg(\alpha\wedge\beta)}$$
From $(\Rightarrow\neg\rightarrow)$, $(\Rightarrow\neg\rightarrow)'$ and Principle 3, we obtain  
$$\mbox{\rm($\Rightarrow\neg\rightarrow)''$ \, } \displaystyle \frac{\Gamma\Rightarrow\Delta, \alpha \hspace{.3cm} \Gamma, \beta\Rightarrow\Delta, \neg\alpha\,\,\,\,\,\, \Gamma, \beta\Rightarrow\Delta, \neg\beta} {\Gamma\Rightarrow\Delta, \neg(\alpha\rightarrow\beta)}$$
Finally, from $(\Rightarrow\neg)$, $(\Rightarrow\neg)'$ and Principle 3,  we get  $\mbox{\rm($\Rightarrow\neg)''$ \, } \displaystyle \frac{\Gamma, \alpha\Rightarrow\Delta, \neg\alpha} {\Gamma\Rightarrow\Delta, \neg\alpha}$.

\begin{defi} Let $\gciore '$ the sequent calculus obtained from \gciore\ by performing the following replacements.
\begin{itemize}
\item $ \mbox{\rm($\Rightarrow\neg\vee)$}$ \, by \, $\mbox{\rm($\Rightarrow\neg\vee)''$}$.
\item $ \mbox{\rm($\Rightarrow\neg\wedge)$}$ \, by \, $\mbox{\rm($\Rightarrow\neg\wedge)''$}$.
\item $ \mbox{\rm($\Rightarrow\neg\to)$}$ \, by \, $\mbox{\rm($\Rightarrow\neg\to)''$}$.
\item $ \mbox{\rm($\Rightarrow\neg)$}$ \, by \, $\mbox{\rm($\Rightarrow\neg)''$}$.
\end{itemize}
\end{defi}

\begin{table}
\centering
{\bf Axioms}
$$\alpha \Rightarrow \alpha$$
{\bf Structural Rules}
$$\mbox{\rm(w$\Rightarrow$) \, }  \displaystyle \frac{\Gamma\Rightarrow\Delta} {\Gamma,\alpha\Rightarrow\Delta} \hspace{2cm}  \mbox{\rm($\Rightarrow$w) }  \displaystyle \frac{\Gamma\Rightarrow\Delta} {\Gamma\Rightarrow\Delta, \alpha}$$
{\bf Logical Rules}
$$ \mbox{\rm($\vee\Rightarrow)$ \, } 
   \displaystyle \frac{\Gamma,\alpha\Rightarrow\Delta \hspace{.5cm} \Gamma,\beta\Rightarrow\Delta} {\Gamma,\alpha\vee\beta\Rightarrow\Delta} \hspace{1.5cm} \mbox{\rm($\Rightarrow\vee)$ \, } \displaystyle \frac{\Gamma\Rightarrow\Delta,\alpha, \beta} {\Gamma\Rightarrow\Delta,\alpha\vee\beta}$$
$$ \mbox{\rm($\neg\vee\Rightarrow)$ \, }   \displaystyle \frac{\Gamma,\alpha,\neg\alpha, \beta, \neg\beta \Rightarrow\Delta \hspace{0.5cm} \Gamma,\neg\alpha, \neg\beta \Rightarrow\Delta, \alpha, \beta} {\Gamma, \neg(\alpha\vee\beta)\Rightarrow\Delta}$$
$$ \mbox{\rm($\Rightarrow\neg\vee)$ \, } \displaystyle \frac{\Gamma,\alpha\Rightarrow\Delta, \neg\alpha  \hspace{0.5cm} \Gamma,\alpha \Rightarrow\Delta, \beta \hspace{0.5cm} \Gamma,\alpha \Rightarrow\Delta, \neg\beta \hspace{0.5cm} \Gamma,\beta \Rightarrow\Delta, \alpha \hspace{0.5cm} \Gamma,\beta \Rightarrow\Delta, \neg\alpha \hspace{0.5cm} \Gamma,\beta \Rightarrow\Delta, \neg\beta} {\Gamma\Rightarrow\Delta, \neg(\alpha\vee\beta)}$$

$$ \mbox{\rm($\wedge\Rightarrow)$ \, } 
   \displaystyle \frac{\Gamma,\alpha,\beta\Rightarrow\Delta} {\Gamma,\alpha\wedge\beta\Rightarrow\Delta} \hspace{1.5cm} \mbox{\rm($\Rightarrow\wedge)$ \, } \displaystyle \frac{\Gamma\Rightarrow\Delta,\alpha \hspace{0.5cm} \Gamma\Rightarrow\Delta,\beta} {\Gamma\Rightarrow\Delta,\alpha\wedge\beta}$$
$$ \mbox{\rm($\neg\wedge\Rightarrow)$ \, } 
  \displaystyle \frac{\Gamma\Rightarrow\Delta, \alpha, \beta \hspace{0.5cm} \Gamma,\neg\alpha \Rightarrow\Delta, \alpha \hspace{0.5cm} \Gamma,\neg\beta \Rightarrow\Delta, \beta \hspace{0.5cm} \Gamma,\neg\alpha, \neg\beta \Rightarrow\Delta} {\Gamma, \neg(\alpha\wedge\beta)\Rightarrow\Delta} $$
 $$ \mbox{\rm($\Rightarrow\neg\wedge)$ \, } \displaystyle \frac{\Gamma,\alpha, \beta\Rightarrow\Delta, \neg\alpha  \hspace{0.5cm} \Gamma,\alpha, \beta \Rightarrow\Delta, \neg\beta} {\Gamma\Rightarrow\Delta, \neg(\alpha\wedge\beta)}$$

$$ \mbox{\rm($\rightarrow\Rightarrow)$ \, } 
   \displaystyle \frac{\Gamma\Rightarrow\Delta,\alpha \hspace{0.5cm} \Gamma,\beta\Rightarrow\Delta} {\Gamma,\alpha\rightarrow\beta\Rightarrow\Delta} \hspace{1.5cm} \mbox{\rm($\Rightarrow\rightarrow)$ \, } \displaystyle \frac{\Gamma,\alpha \Rightarrow\Delta,\beta} {\Gamma\Rightarrow\Delta,\alpha\rightarrow\beta}$$
$$ \mbox{\rm($\neg\rightarrow\Rightarrow)$ \, }     \displaystyle \frac{\Gamma,\alpha,\neg\beta\Rightarrow\Delta, \beta \hspace{0.5cm} \Gamma,\alpha, \neg\alpha, \neg\beta \Rightarrow\Delta} {\Gamma, \neg(\alpha\rightarrow\beta)\Rightarrow\Delta} $$
 $$ \mbox{\rm($\Rightarrow\neg\rightarrow)$ \, } \displaystyle \frac{\Gamma\Rightarrow\Delta, \alpha \hspace{0.5cm} \Gamma, \beta\Rightarrow\Delta, \neg\alpha \hspace{0.5cm} \Gamma, \beta\Rightarrow\Delta, \neg\beta} {\Gamma\Rightarrow\Delta, \neg(\alpha\rightarrow\beta)}$$

$$ \mbox{\rm($\Rightarrow\neg)$ \, } \displaystyle \frac{\Gamma, \alpha\Rightarrow\Delta, \neg\alpha} {\Gamma\Rightarrow\Delta, \neg\alpha}\hspace{1.3cm} \mbox{\rm($\neg\neg\Rightarrow)$ \, } \displaystyle \frac{\Gamma,\alpha\Rightarrow\Delta} {\Gamma,\neg\neg\alpha\Rightarrow\Delta} \hspace{1.3cm} \mbox{\rm($\Rightarrow\neg\neg)$ \, }   \displaystyle \frac{\Gamma\Rightarrow\Delta,\alpha} {\Gamma\Rightarrow\Delta,\neg\neg\alpha}$$
$$ \mbox{\rm($\circ\Rightarrow)$ \, } \displaystyle \frac{\Gamma\Rightarrow\Delta,\alpha \hspace{.5cm} \Gamma\Rightarrow\Delta,\neg\alpha} {\Gamma,\circ\alpha\Rightarrow\Delta} \hspace{1.3cm} \mbox{\rm($\Rightarrow\circ)$ \, } \displaystyle \frac{\Gamma,\alpha, \neg\alpha \Rightarrow\Delta} {\Gamma\Rightarrow\Delta, \circ\alpha}\hspace{1.3cm} \mbox{\rm($\neg\circ\Rightarrow)$ \, } \displaystyle \frac{\Gamma, \alpha, \neg\alpha \Rightarrow\Delta} {\Gamma, \neg{\circ}\alpha\Rightarrow\Delta}$$

\caption{The ${\cal G}Ciore'$ calculus}
\label{tab:abc}
\end{table}

\noindent Then, it is immediate the following proposition.

\begin{prop}\label{prop2.4.6}
Let $\Gamma$ and $\Delta$ be finite sets of formulas, the following conditions are equivalent:
\begin{enumerate}[{\rm (i)}]
\item $\vdash_{\gciore}\Gamma\Rightarrow\Delta$
\item $\vdash_{{\gciore}'}\Gamma\Rightarrow\Delta$
\end{enumerate}
Moreover $\Gamma\Rightarrow\Delta$ is provable in \gciore\ without using the cut rule iff $\Gamma\Rightarrow\Delta$ is provable in \gciore $'$ without using cut.
\end{prop}
\begin{dem} By Remark \ref{rem1} and construction of \gciore $'$.
\end{dem}

\

\noindent As usual, a rule \, $ \displaystyle \frac{\Gamma_{1}\Rightarrow\Delta_{1} \hdots \Gamma_{n}\Rightarrow\Delta_{n}} {\Gamma\Rightarrow\Delta}$ \, is invertible (in ${\cal M}_{e}$)  
if it is verified that: if $\models_{{\cal M}_{e}}\Gamma\Rightarrow\Delta$, then $\models_{{\cal M}_{e}}\Gamma_{i}\Rightarrow\Delta_{i}$ for all $1\leq i \leq n$.\

\begin{lem} (Inversion Principle) All the logic rules of \gciore $'$ are invertible.
\end{lem}
\begin{dem} We shall just check it for $\mbox{\rm($\Rightarrow\neg\vee)''$}$ since the proof is similar for the remaining logic rules of \gciore $'$. Suppose that $\models_{{\cal M}_{e}}\Gamma\Rightarrow\Delta, \neg(\alpha\vee\beta)$ and let $v$ be a ${\cal M}_{e}$-valuation. Since $v$ satisfies $\Gamma\Rightarrow\Delta, \neg(\alpha\vee\beta)$, then $v(\neg(\alpha\vee\beta))\in{\cal D}$. We have the next cases: \\
\underline{Case 1}: $v(\alpha\vee\beta)=\0$. By the definition of $\hat{\vee}$, we have that $v(\alpha)=\0$ and $v(\beta)=\0$. From $v(\alpha)=0$, we have $\models_{{\cal M}_{e}}\Gamma,\alpha\Rightarrow\Delta, \neg\alpha$, $\models_{{\cal M}_{e}}\Gamma,\alpha \Rightarrow\Delta, \beta$  and $\models_{{\cal M}_{e}}\Gamma,\alpha \Rightarrow\Delta, \neg\beta$. On the other hand, from $v(\beta)=\0$ we have $\models_{{\cal M}_{e}}\Gamma,\beta \Rightarrow\Delta, \alpha$, $\models_{{\cal M}_{e}}\Gamma,\beta \Rightarrow\Delta, \neg\alpha$  and $\models_{{\cal M}_{e}}\Gamma,\beta \Rightarrow\Delta, \neg\beta$.\\
\underline{Case 2}: $v(\alpha\vee\beta)=\um$. Then,  $v(\alpha)=\um$ and $v(\beta)=\um$. From $v(\alpha)=\um$ we have $\models_{{\cal M}_{e}}\Gamma,\beta \Rightarrow\Delta, \alpha$. And from $v(\neg\alpha)=\um$ we have $\models_{{\cal M}_{e}}\Gamma,\beta \Rightarrow\Delta, \neg\alpha$ and $\models_{{\cal M}_{e}}\Gamma,\alpha\Rightarrow\Delta, \neg\alpha$. On the other hand, from $v(\beta)=\um$ we have $\models_{{\cal M}_{e}}\Gamma,\alpha \Rightarrow\Delta, \beta$. Since $v(\neg\beta)=\um$ we have $\models_{{\cal M}_{e}}\Gamma,\alpha \Rightarrow\Delta, \neg\beta$  and  $\models_{{\cal M}_{e}}\Gamma,\beta \Rightarrow\Delta, \neg\beta$.
\end{dem}

\

\noindent Recall that a \textit{literal} is a formula that is a propositional variable or a negated propositional variable, i.e., $\alpha$ is a literal if there is a propositional variable $p$ such that $\alpha$ is $p$ or $\alpha$ is $\neg p$. \\
Now, we assign a non negative integer to each formula $\alpha$ of the language of  \ciore\ in the following way.

\begin{defi} (Weight of a formula) Let $\alpha$ be a formula, $w(\alpha)$ is the non negative integer obtained as follows:
\begin{enumerate}[{\rm (i)}]
\item  If $\alpha$ is a literal, $w(\alpha)=0$
\item If $\alpha$ is $\beta\sharp\gamma$, for $\sharp=\{\vee, \wedge, \rightarrow\}$, $w(\alpha)= w(\beta\sharp\gamma)= w(\beta) + w(\gamma) + 1$.
\item If $\alpha$ is $\circ\beta$, $w(\alpha)= w(\circ\beta)= w(\beta) + w(\neg\beta) + 1$.
\item If $\alpha$ is $\neg\neg\beta$, $w(\alpha)= w(\neg\neg\beta)= w(\neg\beta) + 1$.
\item If $\alpha$ is $\neg\circ\beta$, $w(\alpha)= w(\neg\circ\beta)= w(\circ\beta) + 1$.
\item If $\alpha$ is $\neg(\beta\sharp\gamma)$, for $\sharp\in \{\vee, \wedge, \rightarrow\}$, $w(\alpha)= w(\neg(\beta\sharp\gamma))= w(\beta) + w(\neg\beta) + w(\gamma) + w(\neg\gamma) + 2$.
\end{enumerate}
\end{defi}

\

\noindent It is not difficult to check that the notion of weight is well-defined. 
If $\Gamma\Rightarrow\Delta$ is a (finite) sequent, we will call weight of $\Gamma\Rightarrow\Delta$ (denoted by $w(\Gamma\Rightarrow\Delta)$) to
$$w(\Gamma\Rightarrow\Delta)=\sum \limits_{\gamma\in\Gamma} w(\gamma) + \sum \limits_{\delta\in\Delta} w(\delta).$$

\

\begin{theo}\label{2.4.8}{\rm ({\bf Completeness of \gciore $'$ with respect to ${\cal M}_{e}$})} \\ If \, $\models_{{\cal M}_{e}}\Gamma\Rightarrow\Delta$ \, then there exists a cut-free proof of the sequent \, $\Gamma\Rightarrow\Delta$ \, in \gciore $'$.
\end{theo}
\begin{dem}
Let $\Gamma\Rightarrow\Delta$ be valid in ${\cal M}_{e}$. We use induction on the weight of the sequent $\Gamma\Rightarrow\Delta$.\\
\textbf{ Base step:} $w(\Gamma\Rightarrow\Delta)=0$. So all the formulas of $\Gamma\cup\Delta$ are literals. Since the sequent is valid, necessarily $\Gamma\cap\Delta\neq\emptyset$. Otherwise, we can always find a valuation that refutes $\Gamma\Rightarrow\Delta$. Let $p\in \Gamma\cap\Delta$ or $\neg p \in \Gamma\cap\Delta$ then $$ \mbox{\rm(w's) \, } \displaystyle \frac{p \Rightarrow p}{\Gamma\Rightarrow\Delta} \hspace{1.3cm} \mbox{\rm(w's) \, } \displaystyle \frac{\neg p \Rightarrow\neg p}{\Gamma\Rightarrow\Delta},$$ and in both cases $\Gamma\Rightarrow\Delta$ is provable.\\
\textbf{(I.H)} Assume that every valid sequent $\Gamma'\Rightarrow\Delta'$ such that $w(\Gamma'\Rightarrow\Delta')<k$, $k>0$, is provable in \gciore $'$ without using the cut rule.\\
Now, let $\Gamma\Rightarrow\Delta$ be a valid sequent such that $w(\Gamma\Rightarrow\Delta)=k\geq 1$. 
Then, there is a formula $\gamma\in\Gamma$ or ($\gamma\in\Delta$) such that $w(\gamma)\geq 1$. Then, $\gamma$ is of the form of one of the following: (a) $\beta\vee\delta$,  (b) $\beta\wedge\delta$, (c) $\beta\rightarrow\delta$, (d) $\neg(\beta\vee\delta)$, (e) $\neg(\beta\wedge\delta)$, (f)  $\neg(\beta\rightarrow\delta)$, (g)  $\neg\neg\beta$, (h)  $\neg\circ\beta$ or (i)  $\circ\beta$. We analyze just case (f), the rest are analogous. \\
Suppose that $\Gamma\Rightarrow\Delta$ is $\Gamma', \neg(\beta\rightarrow\delta)\Rightarrow\Delta$.
Since $\Gamma', \neg(\beta\rightarrow\delta)\Rightarrow\Delta$ is valid and the rule $ \mbox{\rm($\neg\rightarrow\Rightarrow)$ \, }$ 
$$ \displaystyle \frac{\Gamma',\beta, \neg\delta\Rightarrow\Delta, \delta \hspace{.5cm} \Gamma',\beta, \neg\beta, \neg\delta\Rightarrow\Delta} {\Gamma',\neg(\beta\rightarrow\delta)\Rightarrow\Delta}$$
is invertible, we have that $\Gamma',\beta, \neg\delta\Rightarrow\Delta, \delta$ and $\Gamma',\beta, \neg\beta, \neg\delta\Rightarrow\Delta$ are valid. But $w(\Gamma',\beta, \neg\delta\Rightarrow\Delta, \delta)<w(\Gamma',\neg(\beta\rightarrow\delta)\Rightarrow\Delta)$ and $w(\Gamma',\beta, \neg\beta, \neg\delta\Rightarrow\Delta)<w(\Gamma',\neg(\beta\rightarrow\delta)\Rightarrow\Delta)$, so, by \textbf{(I.H)} $\Gamma',\beta, \neg\delta\Rightarrow\Delta, \delta$ and $\Gamma',\beta, \neg\beta, \neg\delta\Rightarrow\Delta$ are provable in \gciore $'$ without using the cut rule. Then, we can construct the cut-free proof

\vspace{-0.5cm}
\begin{prooftree}
\AxiomC{$\vdots $}
\LeftLabel{}
\UnaryInfC{$\Gamma',\beta, \neg\delta\Rightarrow\Delta, \delta$}
\AxiomC{$\vdots$}
\LeftLabel{}
\UnaryInfC{$\Gamma',\beta, \neg\beta, \neg\delta\Rightarrow\Delta$}
\LeftLabel{\small($\vee\Rightarrow$ )}
\BinaryInfC{$\Gamma',\neg(\beta\rightarrow\delta)\Rightarrow\Delta$}
\end{prooftree}
of the sequent $\Gamma',\neg(\beta\rightarrow\delta)\Rightarrow\Delta$ 
\end{dem}

\

\begin{cor}\label{cor1}
\gciore $'$ admits cut-elimination.
\end{cor}
\begin{dem}
Immediate consequence of the previous theorem.
\end{dem}

\

\begin{cor}\label{corcutelim} \gciore\  admits cut-elimination.
\end{cor}
\begin{dem}
From Corollary \ref{cor1} and Proposition \ref{prop2.4.6}.
\end{dem}

\

\section{Some consequences of cut elimination for \ciore}\label{s5}

In this section, we show some consequences of the cut-elimination theorem for \ciore. In first place, we introduce the notion of generalize subformula in order to show that there exists a syntactic decision procedure for \ciore .

\begin{defi}\label{defgsub}{\rm (Generalized subformula for \ciore )}
The set of generalized subformulas of a given formula $\gamma$, $\bf{gsub}(\gamma)$, is defined as the least set of formulas fulfilling the following conditions:
\begin{enumerate}[\rm (1)]
\item $\alpha\in\bf{gsub}(\alpha)$.
\item $\bf{gsub}(\alpha)\subseteq\bf{gsub}(\neg\alpha)$.
\item $\bf{gsub}(\alpha)\cup\bf{gsub}(\beta)\subseteq\bf{gsub}(\alpha\#\beta)$ where $\#\in\{\wedge , \vee , \rightarrow\}$.
\item $\bf{gsub}(\neg\alpha)\cup\bf{gsub}(\neg\beta)\subseteq\bf{gsub}(\neg(\alpha\#\beta))$ where $\#\in\{\wedge , \vee , \rightarrow\}$.
\item $\bf{gsub}(\neg\alpha)\subseteq\bf{gsub}(\circ\alpha)$.
\end{enumerate}
\end{defi}

\begin{rem}\label{remgsub} From the above definition we may conclude that for every formula $\alpha$ we have that
\begin{enumerate}[{\rm (i)}]
\item $\alpha$ and $\neg\alpha$ are generalized subformulas of $\circ \alpha$,
\item $\alpha$ and $\neg\alpha$ is a generalized subformula $\neg{\circ}\alpha$,
\item $\alpha$ is a generalized subformula of $\neg\neg\alpha$.
\end{enumerate}
\end{rem}

\begin{prop}{\rm (Generalized subformula property)}
Let ${\cal D}$ be a cut-free derivation of $\Gamma\Rightarrow\Delta$ in \gciore. For every sequent $\Pi\Rightarrow\Lambda$ occurring in  ${\cal D}$ it holds:\\
$$\Pi\cup\Lambda \subseteq {\bf gsub}(\Gamma\Rightarrow\Delta)$$
\end{prop}
\begin{dem}The proof is using induction on the number of rule applications in ${\cal D}$, Definition \ref{defgsub}, Remark \ref{remgsub} and inspecting the rules of \gciore . 
\end{dem}

\

\noindent  It is clear that for any sequent $\Gamma\Rightarrow\Delta$, ${\bf gsub}(\Gamma\Rightarrow\Delta)$ is a finite set. Therefore, it is not difficult to state a decision procedure for \gciore\ adapting the procedure given by Gentzen for {\bf PK} (the propositional version of his well-known  {\bf LK}). Therefore,

\begin{theo} \ciore\ has a decision procedure which guarantees bottom-up proof search.
\end{theo}

\noindent Next, we see that, although \ciore\ is a paraconsistent logic, it does not entail contradictions. More precisely, we prove that in \gciore\ no contradiction can be proved.

\begin{lem}\label{lemcon}
Let $\alpha$ be a formula. The following conditions are equivalent:
\begin{enumerate}[\rm (i)]
\item $\Rightarrow\alpha\wedge\neg\alpha$ is provable.
\item The empty sequent \, $\Rightarrow$ \, is provable.
\end{enumerate}
\end{lem}
\begin{dem}
(i)$\Rightarrow$ (ii):  We use induction on the complexity, $c(\alpha)$,  of $\alpha$:\\
If $c(\alpha)=0$, $\alpha$ is a propositional variable. Suppose that there exists a cut-free proof $\D$ of $\Rightarrow p\wedge\neg p$ as follows

\begin{prooftree}
\AxiomC{$\vdots $}
\LeftLabel{\small($\Rightarrow\wedge)$}
\UnaryInfC{$p\wedge\neg p$}
\end{prooftree}
Then, the last inference occurring in $\D$ has to be $(\Rightarrow\wedge)$ (since $\D$ is a cut-free proof). Therefore, $\D$ has the form
\begin{prooftree}
\AxiomC{$\vdots $}
\LeftLabel{}
\UnaryInfC{$\Rightarrow p$}
\AxiomC{$\vdots$}
\LeftLabel{(r)}
\UnaryInfC{$\Rightarrow\neg p$}
\LeftLabel{}
\BinaryInfC{$\Rightarrow p\wedge\neg p$}
\end{prooftree}
and then, (r) must be $(\Rightarrow\neg)$. So, we have a cut-free proof of $p \Rightarrow$ and we can construct the following proof of the empty sequent. 
\begin{prooftree}
\AxiomC{$\vdots $}
\LeftLabel{}
\UnaryInfC{$\Rightarrow p$}
\AxiomC{$\vdots$}
\LeftLabel{}
\UnaryInfC{$p\Rightarrow$}
\LeftLabel{}
\BinaryInfC{$\Rightarrow$}
\end{prooftree}
{\bf (I.H.)} Assume that if  $\gamma$ is such that $c(\gamma)<k$ then if $\Rightarrow\gamma\wedge\neg\gamma$ is provable, then the empty sequent is also provable.\\
Now, let $\alpha$ be such that $c(\alpha)=k>0$. By hypothesis, there is a cut-free proof $\D$ of $\Rightarrow\alpha\wedge\neg\alpha$. Since $\D$ is cut-free, the last rule must be $(\Rightarrow\wedge)$ therefore $\D$ has the form
\begin{prooftree}
\AxiomC{$\vdots $}
\LeftLabel{}
\UnaryInfC{$\Rightarrow\alpha$}
\AxiomC{$\vdots$}
\LeftLabel{$(r)$}
\UnaryInfC{$\Rightarrow\neg\alpha$}
\LeftLabel{$(\Rightarrow\wedge)$}
\BinaryInfC{$\Rightarrow\alpha\wedge\neg\alpha$}
\end{prooftree}
Then, (r) can be one of the rules:  $(\Rightarrow\neg), (\Rightarrow\neg\neg), (\Rightarrow\neg\vee), (\Rightarrow\neg\wedge)$ \, or \, \break  $(\Rightarrow\neg\rightarrow)$. The case where (r) is $(\Rightarrow\neg)$ is analogous to the basic case.\\
If $(r)$ is $(\Rightarrow\neg\neg)$ then $\D$ is of the form
\begin{prooftree}
\AxiomC{$\vdots $}
\LeftLabel{}
\UnaryInfC{$\Rightarrow\neg\beta$}

\AxiomC{$\vdots$}
\UnaryInfC{$\Rightarrow\beta$}
\LeftLabel{$(\Rightarrow\neg\neg)$}
\UnaryInfC{$\Rightarrow\neg\neg\beta$}

\LeftLabel{$(\Rightarrow\wedge)$}
\BinaryInfC{$\Rightarrow\neg\beta\wedge\neg\neg\beta$}
\end{prooftree}
Then, we have (cut-free) proofs of $\Rightarrow\beta$ and $\Rightarrow\neg\beta$ and using $(\Rightarrow\wedge)$ we have a proof for the sequent $\Rightarrow\beta\wedge\neg\beta$. Since $c(\beta)<c(\alpha)=k$, by {\bf (H.I.)} we have that \, $\Rightarrow$ \, is provable. \\
If $(r)$ is $(\Rightarrow\neg\vee)$, $\D$ is of the form:
\begin{prooftree}
\AxiomC{$\vdots $}
\LeftLabel{}
\UnaryInfC{$\Rightarrow\alpha\vee\beta$}

\AxiomC{$\vdots$}
\UnaryInfC{$\Rightarrow\alpha$}
\AxiomC{$\vdots$}
\UnaryInfC{$\Rightarrow\neg\alpha$}
\AxiomC{$\vdots$}
\UnaryInfC{$\Rightarrow\beta$}
\AxiomC{$\vdots$}
\UnaryInfC{$\Rightarrow\neg\beta$}
\LeftLabel{$(\Rightarrow\neg\vee)$}
\QuaternaryInfC{$\Rightarrow\neg(\alpha\vee\beta)$}

\LeftLabel{$(\Rightarrow\wedge)$}
\BinaryInfC{$\Rightarrow(\alpha\vee\beta)\wedge\neg(\alpha\vee\beta)$}
\end{prooftree}
Then we have cut-free proofs of $\Rightarrow\beta$ and $\Rightarrow\neg\beta$ and using  $(\Rightarrow\wedge)$ we have that $\Rightarrow\beta\wedge\neg\beta$ is provable. Since $c(\beta)<c(\alpha)=k$, by {\bf (I.H.)}, \, $\Rightarrow$ \, is provable. The cases where $(r)$ is $(\Rightarrow\neg\wedge)$ or $(\Rightarrow\neg\to)$ are analogous. \\[2mm]
(ii)$\Rightarrow$(i): Immediate.
\end{dem}

\

\begin{cor} \ciore\ does not entail contradictions. 
\end{cor}
\begin{dem}
Suppose that there exists $\alpha$ such that $\Rightarrow\alpha\wedge\neg\alpha$ is provable. Then by Lemma \ref{lemcon}, the empty sequent is provable. By Corollary \ref{corcutelim}, there exists a cut-free proof of the empty sequent; and this is impossible.
\end{dem}

\begin{cor} For every theorem $\alpha$ of \ciore , and every $\mm M_{e}$-valuation $v$, $v(\alpha)=\1$.
\end{cor}

\section{First-order case}\label{s6}

In this section, we present a sound a complete sequent system for the first-order version of \ciore, \qciore. Then, we show that this system enjoys the cut-elimination property.

\begin{defi} Let \gqciore\ the sequent system formed by the first-order versions of the rules and axioms of \gciore\ plus the following rules:

\

$$\mbox{\rm ($\forall\Rightarrow$)} \,  \displaystyle \frac{\phi(b),\Gamma\Rightarrow\Delta} {\forall x \phi(x), \Gamma\Rightarrow\Delta}\hspace{2.5cm} \mbox{\rm ($\Rightarrow\forall$)}\, \displaystyle \frac{\Gamma\Rightarrow\Delta, \phi(a)}{\Gamma\Rightarrow\Delta, \forall x \phi(x)}$$

\

\noindent where $b$ is an arbitrary free variable and the eigenvariable $a$ does not occur in the lower sequent. Besides, in ($\Rightarrow\pt$) all occurrences of $a$ in $\phi(a)$ are indicated. 

$$\mbox{\rm ($\exists\Rightarrow$)} \, \displaystyle \frac{\phi(a),\Gamma\Rightarrow\Delta} {\exists x \phi(x), \Gamma\Rightarrow\Delta}\hspace{2.5cm}\mbox{\rm ($\Rightarrow\exists$)} \, \displaystyle \frac{\Gamma\Rightarrow\Delta, \phi(b)} {\Gamma\Rightarrow\Delta, \exists x \phi(x)}$$

\

\noindent where the eigenvariable $a$ does not occur in the lower sequent and $b$ is an arbitrary free variable.

$$\mbox{\rm ($\circ\exists\Rightarrow$)}\, \displaystyle \displaystyle \frac{\circ\phi(a),\Gamma\Rightarrow\Delta} {\circ\exists x \phi(x), \Gamma\Rightarrow\Delta} \hspace{2.5cm}  \mbox{\rm ($\Rightarrow\circ\exists$)}\, \displaystyle \frac{\Gamma\Rightarrow\Delta, \circ\phi(b)} {\Gamma\Rightarrow\Delta, \circ\exists x \phi(x)}$$

\

\noindent where the eigenvariable $a$ does not occur in the lower sequent and $b$ is an arbitrary free variable.

$$\mbox{\rm ($\circ\forall\Rightarrow$)}\, \displaystyle \frac{\circ\phi(b),\Gamma\Rightarrow\Delta} {\circ\forall x \phi(x), \Gamma\Rightarrow\Delta}\hspace{2.5cm}\mbox{\rm ($\Rightarrow\circ\forall$)} \,\displaystyle \frac{\Gamma\Rightarrow\Delta, \circ\phi(b)} {\Gamma\Rightarrow\Delta, \circ\forall x \phi(x)}$$

\noindent where $b$ is an arbitrary free variable.
\end{defi}

\

\begin{prop} The following sequents are provable in \gqciore .
\begin{enumerate}[\rm (i)]
\item $\Rightarrow\phi(b)\rightarrow\exists x\phi(x)$,
\item $\Rightarrow\forall x \phi(x)\rightarrow\phi(b)$,
\item $\Rightarrow\circ\exists x \phi(x)\leftrightarrow\exists x \circ\phi(x)$, 
\item $\Rightarrow\circ\forall x\phi(x)\leftrightarrow\exists x \circ\phi(x)$.
\end{enumerate}
\end{prop}
\begin{dem} (i) 
\begin{prooftree}
\AxiomC{$\phi(b)\Rightarrow\phi(b)$}
\LeftLabel{\small($\Rightarrow\exists$)}
\UnaryInfC{$\phi(b)\Rightarrow\exists x\phi(b)$}
\LeftLabel{\small( $\Rightarrow\rightarrow$ )}
\UnaryInfC{$\Rightarrow\phi(b)\rightarrow\exists x\phi(x)$}
\end{prooftree}

\noindent (ii)

\begin{prooftree}
\AxiomC{$\phi(b)\Rightarrow\phi(b)$}
\LeftLabel{\small($\forall\Rightarrow$)}
\UnaryInfC{$\forall x\phi(x)\Rightarrow\phi(b)$}
\LeftLabel{\small( $\Rightarrow\rightarrow$ )}
\UnaryInfC{$\Rightarrow\forall x \phi(x)\rightarrow\phi(b)$}
\end{prooftree}

\noindent (iii)

\begin{center}

\hspace{-.7cm}
\AxiomC{$\circ\phi(a)\Rightarrow\circ\phi(a)$}
\LeftLabel{\small($\Rightarrow\exists$)}
\UnaryInfC{$\circ\phi(a)\Rightarrow\exists x\circ\phi(x)$}
\LeftLabel{\small( $\circ\exists\Rightarrow$ )}
\UnaryInfC{$\circ\exists x \phi(x)\Rightarrow\exists x \circ\phi(x)$}
\LeftLabel{\small( $\Rightarrow\rightarrow$ )}
\UnaryInfC{$\Rightarrow\circ\exists x \phi(x)\rightarrow\exists x \circ\phi(x)$}
\DisplayProof \hspace{1cm}
\AxiomC{$\circ\phi(a)\Rightarrow\circ\phi(a)$}
\LeftLabel{\small($\Rightarrow\circ\exists$)}
\UnaryInfC{$\circ\phi(a)\Rightarrow\circ\exists x\phi(x)$}
\LeftLabel{\small( $\exists\Rightarrow$ )}
\UnaryInfC{$\exists x \circ\phi(x)\Rightarrow\circ\exists x \phi(x)$}
\LeftLabel{\small( $\Rightarrow\rightarrow$ )}
\UnaryInfC{$\Rightarrow\exists x \circ\phi(x)\rightarrow\circ\exists x \phi(x)$}
\DisplayProof

\end{center}

\noindent (iv)

\begin{center}

\hspace{-.7cm}
\AxiomC{$\circ\phi(b)\Rightarrow\circ\phi(b)$}
\LeftLabel{\small($\circ\forall\Rightarrow$)}
\UnaryInfC{$\circ\forall x\phi(x)\Rightarrow\circ\phi(b)$}
\LeftLabel{\small( $\Rightarrow\exists$ )}
\UnaryInfC{$\circ\forall x \phi(x)\Rightarrow\exists x \circ\phi(x)$}
\LeftLabel{\small( $\Rightarrow\rightarrow$ )}
\UnaryInfC{$\Rightarrow\circ\forall x\phi(x)\rightarrow\exists x \circ\phi(x)$}
\DisplayProof \hspace{1cm}
\AxiomC{$\circ\phi(b)\Rightarrow\circ\phi(b)$}
\LeftLabel{\small($\Rightarrow\circ\forall$)}
\UnaryInfC{$\circ\phi(b)\Rightarrow\circ\forall x\phi(x)$}
\LeftLabel{\small( $\exists\Rightarrow$ )}
\UnaryInfC{$\exists x \circ\phi(x)\Rightarrow\circ\forall x \phi(x)$}
\LeftLabel{\small( $\Rightarrow\rightarrow$ )}
\UnaryInfC{$\Rightarrow\exists x \circ\phi(x)\rightarrow\circ\forall x\phi(x)$}
\DisplayProof 

\end{center}
\end{dem}

\

\begin{prop} The following rules are derivable in \gqciore .
\begin{enumerate}[\rm (i)]
\item $\displaystyle \frac{\Rightarrow\phi\rightarrow\psi(a)}{\Rightarrow\phi\rightarrow\forall x\psi(x)}$ \, where $a$ does not occur in $\phi$, 
\item  $\displaystyle \frac{\Rightarrow \phi(a)\rightarrow\psi}{\Rightarrow \exists x \phi(x)\rightarrow\psi}$ \, where $a$ does not occur in $\psi$.
\end{enumerate}
\end{prop}
\begin{dem}

\noindent (i)

\begin{prooftree}
\AxiomC{}
\LeftLabel{\small(hyp)}
\UnaryInfC{$\Rightarrow\phi\rightarrow\psi(a)$}
\LeftLabel{\small(w $\Rightarrow$)}
\UnaryInfC{$\phi\Rightarrow\phi\rightarrow\psi(a)$}
\LeftLabel{\small( $\Rightarrow$ w)}
\UnaryInfC{$\phi\Rightarrow\phi\rightarrow\psi(a), \psi(a)$}

\AxiomC{$\phi\Rightarrow\phi$}
\LeftLabel{\small($\Rightarrow$ w)}
\UnaryInfC{$\phi\Rightarrow\phi, \psi(a)$}

\AxiomC{$\psi(a)\Rightarrow\psi(a)$}
\LeftLabel{\small(w $\Rightarrow$)}
\UnaryInfC{$\phi , \psi(a)\Rightarrow\psi(a)$}

\LeftLabel{\small( $\rightarrow\Rightarrow$ )}
\BinaryInfC{$\phi , \phi\rightarrow\psi(a) \Rightarrow\psi(a)  $}

\LeftLabel{\small( $cut$ )}
\BinaryInfC{$\phi \Rightarrow\psi(a) $}
\LeftLabel{\small( $\Rightarrow\forall$ )}
\UnaryInfC{$\phi \Rightarrow\forall x\psi(x) $}
\LeftLabel{\small( $\Rightarrow \rightarrow$ )}
\UnaryInfC{$\Rightarrow\phi\rightarrow\forall x\psi(x) $}

\end{prooftree}

\noindent (ii)

\begin{prooftree}
\AxiomC{}
\LeftLabel{\small(hyp)}
\UnaryInfC{$\Rightarrow\phi(a)\rightarrow\psi$}
\LeftLabel{\small(w $\Rightarrow$)}
\UnaryInfC{$\phi(a)\Rightarrow\phi(a)\rightarrow\psi$}
\LeftLabel{\small( $\Rightarrow$ w)}
\UnaryInfC{$\phi(a)\Rightarrow\phi(a)\rightarrow\psi, \psi$}

\AxiomC{$\phi(a)\Rightarrow\phi(a)$}
\LeftLabel{\small($\Rightarrow$ w)}
\UnaryInfC{$\phi(a)\Rightarrow\phi(a), \psi$}

\AxiomC{$\psi\Rightarrow\psi$}
\LeftLabel{\small(w $\Rightarrow$)}
\UnaryInfC{$\phi(a) , \psi\Rightarrow\psi$}

\LeftLabel{\small( $\rightarrow\Rightarrow$ )}
\BinaryInfC{$\phi(a) , \phi(a)\rightarrow\psi \Rightarrow\psi  $}

\LeftLabel{\small( $cut$ )}
\BinaryInfC{$\phi(a) \Rightarrow\psi $}
\LeftLabel{\small( $\exists\Rightarrow$ )}
\UnaryInfC{$\exists x\phi(x) \Rightarrow\psi $}
\LeftLabel{\small( $\Rightarrow \rightarrow$ )}
\UnaryInfC{$\Rightarrow\exists x\phi(x)\rightarrow\psi $}

\end{prooftree}

\end{dem}

\

\begin{defi}\label{defSeqValid}
Let $\langle\ff A,\tripc{\cdot}\rangle$ be a \qciore-structure, and let $\Gamma\Rightarrow\Delta$ be a \gqciore -sequent. We say that $\Gamma\Rightarrow\Delta$ is satisfied in $\ff A$ by the assignment $s$ if either some formula in $\Gamma$ is not satisfied by $s$ (in $\ff A$), or some formula in $\Delta$  is satisfied by $s$ (in $\ff A$). A sequent is valid in  $\ff A$ if it is satisfied by every assignment $s$ in $\ff A$. Finally, a sequent $\Gamma\Rightarrow\Delta$ is valid if it is satisfied by all \qciore-structures and, in this case, we write $\models\Gamma\Rightarrow\Delta$.
\end{defi}

\begin{theo}\label{soundQ} (Soundness) Let $\Gamma, \Delta$ be a set of formulas.
$$ \mbox{If } \Gamma\Rightarrow\Delta \mbox{ is provable \gqciore, then } \models_\qciore \Gamma\Rightarrow\Delta.$$
\end{theo}
\begin{dem}  As usual, the proof is by induction on the number of instances of rule applications in the derivation of $\Gamma\Rightarrow\Delta$. So, we just have to check that all the rules of \gqciore\ preserve validity. Next, we check it just for ($\neg\vee\Rightarrow$), ($\circ\forall\Rightarrow$) and $(\Rightarrow\circ\exists)$; the rest of the rules of are analyzed similarly.
Let us prove  that 
$$(\neg\vee\Rightarrow)  \displaystyle \frac{\Gamma,\alpha,\neg\alpha, \beta, \neg\beta \Rightarrow\Delta \hspace{0.5cm} \Gamma,\neg\alpha, \neg\beta \Rightarrow\Delta, \alpha, \beta} {\Gamma, \neg(\alpha\vee\beta)\Rightarrow\Delta}$$
 preserves validity. Suppose that $\Gamma,\alpha,\neg\alpha, \beta, \neg\beta \Rightarrow\Delta$ and $\Gamma,\neg\alpha, \neg\beta \Rightarrow\Delta, \alpha, \beta$  are valid. Then, for every structure $\ff A$ and every assignment $s$, $s$ satisfies $\Gamma,\alpha,\neg\alpha, \beta, \neg\beta \Rightarrow\Delta$ and $\Gamma,\neg\alpha, \neg\beta \Rightarrow\Delta, \alpha, \beta$. If either $s$ does not satisfy some formula of $\Gamma$ or $s$ satisfies all formulas of $\Delta$, the proof is completed. Otherwise, $s$ satisfies all formulas of $\Gamma$ and none of $\Delta$. Let us see that, in this case, $s$ does not satisfies $\neg(\alpha\vee\beta)$. Suppose that $s$ satisfies $\neg(\alpha\vee\beta)$. Then $s\in \mas{\neg(\alpha\vee\beta)}\cup\por{\neg(\alpha\vee\beta)}$ and we have the following alternatives:
\begin{itemize}
\item If $s \in \mas{\neg(\alpha\vee\beta)}$ then we have that $s \in  \menos{\alpha\vee\beta}$ and then $s \in \menos\alpha\cap\menos\beta$. So, $s \in \menos\alpha$ and $s\in\menos\beta$, i.e.  $s$ does not satisfy $\alpha$ and $s$ does not satisfy $\beta$, which is a contradiction since $s$ satisfies $\Gamma,\neg\alpha, \neg\beta \Rightarrow\Delta, \alpha, \beta$. 
\item  If $s \in \por{\neg(\alpha\vee\beta)}$ then $s \in  \por{\alpha\vee\beta}$, i.e. $s \in \por\alpha\cap\por\beta$. So, $s \in \por\alpha$ and $s\in\por\beta$ and we have that $s$ satisfies $\alpha$ and $\neg\alpha$; and $s$ satisfies $\beta$ and $\neg\beta$, which is a contradiction since $s$ satisfies $\Gamma,\alpha,\neg\alpha, \beta, \neg\beta \Rightarrow\Delta$.
\end{itemize}
Therefore $s$ does not satisfy  $\neg(\alpha\vee\beta)$ and then $s$ satisfies $\Gamma, \neg(\alpha\vee\beta)\Rightarrow\Delta$.

In order to see that  $(\circ\forall\Rightarrow) \displaystyle \frac{\circ\phi(a),\Gamma\Rightarrow\Delta} {\circ\forall x \phi(x), \Gamma\Rightarrow\Delta}$  preserves validity, suppose that ${\circ}\phi(a),\Gamma\Rightarrow\Delta$ is valid. Then, for every structure $\ff A$ and every assignment $s$, $s$ satisfies ${\circ}\phi(a),\Gamma\Rightarrow\Delta$. If $s$ does not satisfy some formula of $\Gamma$ or $s$ satisfies any every formula of $\Delta$, the proof is completed. 
Otherwise, $s$ does not satisfy ${\circ}\phi(a)$.  Then $\menos{{\circ}\phi(a)}=S(\ff A)$, that is $\por{\phi(a)}=S(\ff A)$. If it were the case that $s$ satisfies ${\circ}\forall x \phi(x)$, then $s\in\mas{{\circ}\forall x \phi(x)}\cup\por{{\circ}\forall x \phi(x)}$. Since $\por{{\circ}\psi}=\emptyset$ for all $\psi$, we have that $s\in\mas{{\circ}\forall x \phi(x)}$. That is, $s\in\mas{\forall x \phi(x)}\cup\menos{\forall x \phi(x)}$ and so,
\begin{itemize}
\item if $s\in\mas{\forall x \phi(x)}$ then $s\in \widehat{\ex a}(\mas{\phi(a)})-\widehat{\ex a}(\menos{\phi(a)})$. That is, $s=s_{a}^{m}\in\mas{\phi(a)}$ for some $m\in A (m=s(a))$ which is a contradiction since $\mas{\phi(a)}=\emptyset$;
\item if $s \in \menos{\forall x \phi(x)}$ then $s \in \widehat{\ex a}(\menos{\phi(a)})$ and then  $s=s_{a}^{m}\in\menos{\phi(a)}$ for some $m\in A (m=s(a))$, which is a contradiction since $\menos{\phi}=\emptyset$.
\end{itemize}
Therefore $s$ does not satisfy  ${\circ}\forall x \phi(x)$ and then $s$ satisfies ${\circ}\forall x \phi(x), \Gamma\Rightarrow\Delta$.

To see that $(\Rightarrow\circ\exists) \displaystyle \frac{\Gamma\Rightarrow\Delta, {\circ}\phi(a)} {\Gamma\Rightarrow\Delta, {\circ}\exists x \phi(x), }$  preserves validity, suppose that $\Gamma\Rightarrow\Delta, {\circ}\phi(a)$ is valid. 
Then, for every structure $\ff A$ and every assignment $s$, $s$ satisfies $\Gamma\Rightarrow\Delta, {\circ}\phi(a)$.  If $s$ does not satisfy some formula of $\Gamma$ or $s$ satisfies any every formula of $\Delta$, the proof is completed. Otherwise, $s$ satisfies ${\circ}\phi(a)$. \\
Then $\mas{{\circ}\phi(a)}\cup\por{{\circ}\phi(a)}=S(\ff A)$, since $\por{{\circ}\phi(a)}=\emptyset$, then  $\mas{{\circ}\phi(a)}=S(\ff A)$ that is $\mas{\phi(a)}\cup\menos{\phi(a)}=S(\ff A)$ then $\por{\phi(a)}=\emptyset$. Let us suppose that $s$ does not satisfy $\circ\exists x \phi (x)$ then $s\in\menos{\circ\exists x \phi (x)}$ that is $s\in\por{\exists x \phi (x)}$ then $s\in\widehat{\forall} (a)(\por{\phi(a)})$. That is $s=s_{a}^{m}\in\por{\phi(a)}$ for $m=s(a)$ which is a contradiction since $\por{\phi(a)}=\emptyset$.\end{dem}

\


Next, we prove the completeness and cut-elimination theorems \gqciore\ by  using the well-known Sch\"utte's method (see~\cite{Tau}).

Recall that an expression $\Gamma\Rightarrow\Delta$ is called an {\em infinite sequent} if $\Gamma$ and $\Delta$ are infinite (countable) sets of formulas. An infinite sequent $\Gamma\Rightarrow\Delta$ is called provable if a finite part $\Gamma'\Rightarrow\Delta'$ of the sequent is provable, i.e., $\Gamma'\subseteq \Gamma$ and $\Delta'\subseteq \Delta$ are finite.

We define, for each sequent $S$, a (possible infinite) tree, called the {\em reduction tree} for $S$, from which we can obtain either a cut-free proof of $S$ or a \qciore-structure not satisfying $S$. This method is due to Sch\"utte (see \cite{Tau}). 
This reduction tree for S, denoted by $\T(S)$, contains a sequent at each node.

\begin{defi} \label{defi6.6} {\rm (Reduction tree of $S$)}  It is constructed in stages as follows. \\
{\em Stage 0:} Write $S$ at the bottom of the tree.\\
{\em Stage $k$ $(k>0)$:} This is defined by cases:\\
{\bf Case I:} Every topmost sequent has a formula common to its antecedent and succedent. Then stop.\\
{\bf Case II:} Not case I. Then this stage is defined according as\\
$$k \equiv 0, 1, \hdots , 26 \, \, (\mbox{mod } 27).$$
In order to make it simpler, let us assume that there are no individual or function constants. All the free variables which occur in any sequent which has been obtained at or before stage $k$ are {\bf available} at stage $k$. In case there is none, pick any free variable and said that it is available.\\
\begin{enumerate}[\rm (1)]
\item  $k\equiv 0$. ($\circ\Rightarrow$)-reduction. Let $\Pi\Rightarrow\Lambda$ be any topmost sequent of the tree which has been defined by stage $k-1$. Let $\circ \alpha_{1},\hdots , \circ \alpha_{n}$ be all formulas in $\Pi$ whose outermost logical symbol is $\circ$, and to which no reduction has been applied in previous stages. Then write down all sequents of the form
$$\Pi\Rightarrow\Lambda , \gamma_{1},\hdots, \gamma_{n},$$ 
where $\gamma_i$ is $\alpha_i$ or $\neg \alpha_i$, above $\Pi\Rightarrow\Lambda$. So, there are $2^n$ such sequents above $\Pi\Rightarrow\Lambda$.
\item $k\equiv 1$. ($\Rightarrow\circ$)-reduction. Let $\circ \alpha_{1},\hdots , \circ \alpha_{n}$ be all formulas in $\Lambda$ whose outermost logical symbol is $\circ$, and to which no reduction has been applied in previous stages. Then write down the sequent
$$\Pi , \alpha_{1},\hdots, \alpha_{n}, \neg \alpha_{1},\hdots, \neg \alpha_{n}\Rightarrow\Lambda$$
above $\Pi\Rightarrow\Lambda$. 
\item $k\equiv 2$. ($\Rightarrow\neg$)-reduction. Let $\neg \alpha_{1},\hdots , \neg \alpha_{n}$ be all formulas in $\Lambda$ whose outermost logical symbol is $\neg$, and to which no reduction has been applied yet. Then write down
$$\Pi , \alpha_{1}, \hdots ,\alpha_{n}\Rightarrow\Lambda$$
above $\Pi\Rightarrow\Lambda$. 
\item $k\equiv 3$. ($\neg{\circ}\Rightarrow$)-reduction. Let $\neg{\circ} \alpha_{1},\hdots , \neg{\circ} \alpha_{n}$ be all formulas in $\Pi$ whose outermost logical symbols are $\neg{\circ}$, and to which no reduction has been applied yet. Then write down
$$\Pi , \alpha_{1},\hdots, \alpha_{n}, \neg \alpha_{1},\hdots, \neg \alpha_{n}\Rightarrow\Lambda$$
above $\Pi\Rightarrow\Lambda$. 
\item $k\equiv 4$. ($\wedge\Rightarrow$)-reduction. Let $\alpha_{1}\wedge \beta_{1}, \hdots , \alpha_{n}\wedge \beta_{n}$ be all formulas in $\Pi$ whose outermost logical symbol is $\wedge$, and to which no reduction has been applied in previous stages. Then write down
$$\Pi, \alpha_{1}, \beta_{1},\alpha_{2}, \beta_{2},\hdots , \alpha_{n}, \beta_{n}\Rightarrow \Lambda$$
above $\Pi\Rightarrow\Lambda$.
\item $k\equiv 5$. ($\Rightarrow\wedge$)-reduction. Let $\alpha_{1}\wedge \beta_{1}, \hdots , \alpha_{n}\wedge \beta_{n}$ be all formulas in $\Lambda$ whose outermost logical symbol is $\wedge$, and to which no reduction has been applied in previous stages. Then write down all the sequents of the form
$$\Pi\Rightarrow\Lambda, \gamma_{1}, \hdots , \gamma_{n},$$
where $\gamma_{i}$ is either $\alpha_{i}$ or $\beta_{i}$, above $\Pi\Rightarrow\Lambda$. There are $2^{n}$ such sequents above $\Pi\Rightarrow\Lambda$. 
\item $k\equiv 6$. ($\vee\Rightarrow$)-reduction. This is defined in a manner symmetric to (6).
\item $k\equiv 7$. ($\Rightarrow\vee$)-reduction. This is defined in a manner symmetric to (5).
\item $k\equiv 8$. ($\rightarrow\Rightarrow$)-reduction. Let $\alpha_{1}\rightarrow \beta_{1}, \hdots , \alpha_{n}\rightarrow \beta_{n}$ be all formulas in $\Pi$ whose outermost logical symbol is $\rightarrow$, and to which no reduction has been applied in previous stages. Then write down the following sequents above $\Pi\Rightarrow\Lambda$
$$\Pi, \beta_{i_1},\beta_{i_2}, \dots ,\beta_{i_k}\Rightarrow\Lambda, \alpha_{j_1}, \alpha_{j_2}, \dots , \alpha_{j_{n-k}},$$
where $i_{1}<\hdots<i_{k}$,   $j_{1}<\hdots<j_{n-k}$ and  $(i_{1}, \hdots, i_{k},j_{1},\hdots, j_{n-k})$ is a permutation of $(1, 2, \hdots , n)$. There are $2^{n}$ such sequents above $\Pi\Rightarrow\Lambda$. 
\item $k\equiv 9$. ($\Rightarrow\rightarrow$)-reduction. Let $\alpha_{1}\rightarrow \beta_{1}, \hdots , \alpha_{n}\rightarrow \beta_{n}$ be all formulas in $\Lambda$ whose outermost logical symbol is $\rightarrow$, and to which no reduction has been applied in previous stages. Then write down
$$\Pi, \alpha_{1},\alpha_{2},\hdots,\alpha_{n}\Rightarrow\Lambda, \beta_{1},\beta_{2}, \hdots,\beta_{n}$$
above $\Pi\Rightarrow\Lambda$. 
\item $k\equiv 10$. ($\neg\vee\Rightarrow$)-reduction. Let $\neg(\alpha_{1}\vee \beta_{1}), \hdots, \neg(\alpha_{n}\vee \beta_{n})$ be all formulas in $\Pi$ whose outermost logical symbols are $\neg\vee$, and to which no reduction has been applied in previous stages. Then write down
$$\alpha_{i_{1}},\hdots,\alpha_{i_{k}}, \beta_{i_{1}},\hdots,\beta_{i_{k}},\neg \beta_{i_{1}},\hdots,\neg \beta_{i_{k}},\neg \alpha_{1},\hdots,\neg \alpha_{n}, \Pi\Rightarrow\Lambda, \alpha_{j_{1}},\hdots,\alpha_{j_{n-k}}, \beta_{j_{1}},\hdots, \beta_{j_{n-k}}$$
where $i_{1}<\hdots<i_{k}$,   $j_{1}<\hdots<j_{n-k}$ and  $(i_{1}, \hdots, i_{k},j_{1},\hdots, j_{n-k})$ is a permutation of $(1, 2, \hdots , n)$. There are $2^{n}$ such sequents above $\Pi\Rightarrow\Lambda$. 

\item $k\equiv 11$. ($\Rightarrow\neg\vee$)-reduction. Let $\neg(\alpha_{1}\vee \beta_{1}), \hdots, \neg(\alpha_{n}\vee \beta_{n})$ be all formulas in $\Lambda$ whose outermost logical symbol is $\neg\vee$, and to which no reduction has been applied in previous stages. Then write down all the sequents 
$$\Pi\Rightarrow\Lambda, \gamma_{1},\hdots, \gamma_{n}$$
where $\gamma_{i}$ is one of the formulas $\alpha_{i}$, $\neg  \alpha_{i}$, $\beta_{i}$ or $\neg \beta_{i}$ ( $1\leq i \leq n$). Taking all possible combinations of such we have $4^n$ such sequents above $\Pi\Rightarrow\Lambda$.
\item $k\equiv 12$.  ($\neg\wedge\Rightarrow$)-reduction. Let $\neg(\alpha_{1}\wedge \beta_{1}), \hdots, \neg(\alpha_{n}\wedge \beta_{n})$ be all formulas in $\Pi$ whose outermost logical symbols are $\neg\wedge$, and to which no reduction has been applied in previous stages. Then write down all the sequents
$$S', \Pi\Rightarrow\Lambda, S$$
above $\Pi\Rightarrow\Lambda$, where, given the set $F_{AB}=\{\alpha_{1},\hdots, \alpha_{n},\beta_{1},\hdots, \beta_{n}\}$, $S\subseteq F_{AB}$, and $S'=\{\gamma_{i}: \gamma_{i}=\neg \alpha_{i}$ if $ \beta_{i}\notin S$ or $\gamma_{i}=\neg \beta_{i}$ if $\alpha_{i}\notin S\}$. It is clear that $|F_{AB}|=2n$, and therefore, there are \break $|\mathcal{P}(F_{AB})|=2^{|F_{AB}|}=2^{2n}=4^{n}$ such sequents above $\Pi\Rightarrow\Lambda$.

\item $k\equiv 13$. ($\Rightarrow\neg\wedge$)-reduction.  Let $\neg(\alpha_{1}\wedge \beta_{1}), \hdots, \neg(\alpha_{n}\wedge \beta_{n})$ be all formulas in $\Lambda$ whose outermost logical symbol is $\neg\wedge$, and to which no ($\Rightarrow\neg\wedge$)-reduction has been applied in previous stages. Then write down all the sequents 
$$\Pi,\alpha_{1}, \dots, \alpha_{n}, \beta_{1}, \dots, \beta_{n}\Rightarrow\Lambda, \gamma_{1},\hdots, \gamma_{n}$$
where $\gamma_{i}$ is one of the formulas $\neg \alpha_{i}$, or $\neg \beta_{i}$ ( $1\leq i \leq n$). Taking all possible combinations of such we have $2^n$ such sequents above $\Pi\Rightarrow\Lambda$.

\item $k\equiv 14$. ($\neg\rightarrow\Rightarrow$)-reduction. Let $\neg(\alpha_{1}\rightarrow \beta_{1}), \hdots, \neg(\alpha_{n}\rightarrow \beta_{n})$ be all formulas in $\Pi$ whose outermost logical symbols are $\neg\rightarrow$, and to which no reduction has been applied in previous stages. Then write down the sequents 
$$\alpha_{1},\hdots , \alpha_{n},\neg \beta_{1},\hdots, \neg \beta_{n}, \neg \alpha_{i_{1}},\hdots,\neg \alpha_{i_{k}}, \Pi \Rightarrow \Lambda, \beta_{j_{1}},\hdots,\beta_{j_{n-k}}$$
where $i_{1}<\hdots<i_{k}$,   $j_{1}<\hdots<j_{n-k}$ and  $(i_{1}, \hdots, i_{k},j_{1},\hdots, j_{n-k})$ is a permutation of $(1, 2, \hdots , n)$. There are $2^{n}$ such sequents above $\Pi\Rightarrow\Lambda$.

\item $k\equiv 15$. ($\Rightarrow\neg\rightarrow$)-reduction.  Let $\neg(\alpha_{1}\rightarrow \beta_{1}), \hdots, \neg(\alpha_{n}\rightarrow \beta_{n})$ be all formulas in $\Lambda$ whose outermost logical symbol is $\neg\rightarrow$, and to which no ($\Rightarrow\neg\rightarrow$)-reduction has been applied in previous stages. Then write down all the sequents 
$$S_{\gamma_1,\dots, \gamma_n}, \Pi \dots, \beta_{n}\Rightarrow\Lambda, \gamma_{1},\hdots, \gamma_{n}$$
where $\gamma_{i}$ is one of the formulas $\alpha_i$, $\neg \alpha_{i}$, or $\neg \beta_{i}$ ( $1\leq i \leq n$); $S_{\gamma_1,\dots, \gamma_n} \subseteq \{\beta_1, \dots, \beta_n\}$ is such that if $\gamma_i=\alpha_i$, then $\beta_i\not\in S_{\gamma_1,\dots, \gamma_n}$. We have $3^n$ such sequents above $\Pi\Rightarrow\Lambda$.

\item $k\equiv 16$. ($\neg\neg\Rightarrow$)-reduction. Let $\neg\neg \alpha_{1}, \hdots, \neg\neg \alpha_{n}$ be all formulas in $\Pi$ whose outermost logical symbols are $\neg\neg$, and to which no reduction has been applied in previous stages. Then write down
$$\alpha_{1},\hdots, \alpha_{n},\Pi\Rightarrow\Lambda$$
above $\Pi\Rightarrow\Lambda$. 

\item $k\equiv 17$. ($\Rightarrow\neg\neg$)-reduction. Let $\neg\neg \alpha_{1}, \hdots, \neg\neg \alpha_{n}$ be all formulas in $\Lambda$ whose outermost logical symbols are $\neg\neg$, and to which no reduction has been applied in previous stages. Then write down
$$\Pi\Rightarrow\Lambda, \alpha_{1},\hdots, \alpha_{n}$$
above $\Pi\Rightarrow\Lambda$. 

\item $k\equiv 18$.   ($\forall\Rightarrow$)-reduction. Let $\forall x_{1} \alpha_{1}(x_{1}), \hdots, \forall x_{n} \alpha_{n}(x_{n})$ be all formulas in $\Pi$ whose outermost logical symbol is $\forall$. Let $a_{i}$ be the first variable available at this stage which has not been used for a reduction of $\forall x_{i} \alpha_{i}(x_{i})$, for $1\leq i\leq n$. Then write down
$$ \alpha_{1}(a_{1}),\hdots, \alpha_{n}(a_{n}), \Pi\Rightarrow\Lambda$$
above $\Pi\Rightarrow\Lambda$. 

\item $k\equiv 19$. ($\Rightarrow\forall$)-reduction. Let $\forall x_{1} \alpha_{1}(x_{1}), \hdots, \forall x_{n} \alpha_{n}(x_{n})$ be all formulas in $\Lambda$ whose outermost logical symbol is $\forall$  and to which no reduction has been applied so far. Let $a_{1},\hdots, a_{n}$ be the first $n$ free variables (in the list of variables) which are not available at this stage. Then write down
$$ \Pi\Rightarrow\Lambda, \alpha_{1}(a_{1}),\hdots, \alpha_{n}(a_{n})$$
above $\Pi\Rightarrow\Lambda$. Note that  $a_{1},\hdots, a_{n}$ are new available free variables.

\item $k\equiv 20$. ($\exists\Rightarrow$)-reduction. This is defined in a symmetric manner to (20).

\item $k\equiv 21$. ($\Rightarrow\exists$)-reduction. This is defined in a symmetric manner to (19).

\item $k\equiv 22$. ($\circ\forall\Rightarrow$)-reduction. Analogous to (19).

\item $k\equiv 23$. ($\Rightarrow\circ\forall$)-reduction. Let $\circ\forall x_{1} \alpha_{1}(x_{1}), \hdots, \circ\forall x_{n} \alpha_{n}(x_{n})$ be all formulas in $\Lambda$ whose outermost logical symbol is $\circ\forall$. Let $a_{i}$ be the first variable available at this stage which has not been used for a reduction of $\circ\forall x_{i} \alpha_{i}(x_{i})$ for $1\leq i\leq n$. Then write down
$$ \Pi\Rightarrow\Lambda, \circ \alpha_{1}(a_{1}),\hdots, \circ \alpha_{n}(a_{n})$$
above $\Pi\Rightarrow\Lambda$. 

\item $k\equiv 24$. ($\circ\exists\Rightarrow$)-reduction.  Let ${\circ}\exists x_{1}\alpha_{1}(x_{1}), \hdots, {\circ}\exists x_{n} \alpha_{n}(x_{n})$ be all formulas in $\Pi$ whose outermost logical symbols are ${\circ}\exists$. Let $a_{1},\hdots, a_{n}$ be the first $n$ free variables (in the list of variables) which are not available at this stage. Then write down
$$ {\circ}\alpha_{1}(a_{1}),\hdots, {\circ}\alpha_{n}(a_{n}), \Pi\Rightarrow\Lambda$$
above $\Pi\Rightarrow\Lambda$.

\item $k\equiv 25$. ($\Rightarrow\circ\exists$)-reduction. Let ${\circ}\exists x_{1}\alpha_{1}(x_{1}), \hdots, {\circ}\exists x_{n} \alpha_{n}(x_{n})$ be all formulas in $\Lambda$ whose outermost logical symbols are ${\circ}\exists$. Let $a_{i}$ be the first variable available at this stage which has not been used for a reduction of  ${\circ}\exists x_{1}\alpha_{1}(x_{1}), \hdots, {\circ}\exists x_{n} \alpha_{n}(x_{n})$, $1\leq i\leq n$.  Then write down
$$  \Pi\Rightarrow\Lambda, {\circ}\alpha_{1}(a_{1}),\hdots, {\circ}\alpha_{n}(a_{n})$$
above $\Pi\Rightarrow\Lambda$. 

\item $k\equiv 26$.  If $\Pi$ and $\Lambda$ have any formula in common, write nothing above $\Pi\Rightarrow\Lambda$ (so this remains a topmost sequent). If $\Pi$ and $\Lambda$ have no formula in common and the reduction described in {\rm (1)}-{\rm (26)} are not applicable, write the same sequent $\Pi\Rightarrow\Lambda$ again above it.
\end{enumerate}
So the collection of those sequents which are obtained by the above reduction process, together with the partial order obtained by this process, is the reduction tree (for $S$) and it is denoted by $\T(S)$.
\end{defi}

\

As usual, a (finite or infinite) sequence $S_0, S_1, \dots$ of sequents of $\T(S)$ is a {\em branch} if (1) $S_0$ is $S$; (2) $S_{i+1}$ stands immediately above $S_i$; (3) if the sequence is finite, say $S_1, \dots, S_n$, then there exists at least one formula in common in the antecedent and the succedent of $S_n$.

\begin{lem} Let $S$ be a sequent. Then either there is a cut-free proof of $S$, or there is a \qciore-structure which refutes $S$.
\end{lem}
\begin{dem} Let $\T(S)$ be the reduction tree of $S$. If each branch of $\T(S)$ is finite and ends with a sequent whose antecedent and succedent contain a formula in common, then it is a routine task to write a cut-free proof for $S$. Otherwise, there is an infinite branch of $\T(S)$. By K\"onig's lemma, we have an infinite branch in $\T(S)$, consisting of  
$$\Gamma_0\Rightarrow \Delta_0, \,  \Gamma_1\Rightarrow \Delta_1, \, \dots, \, \Gamma_i\Rightarrow \Delta_i, \, \dots$$
Let \, $\Gamma = \bigcup \limits_{i=0}^{\infty} \Gamma_i$ \, and \, $\Delta = \bigcup \limits_{i=0}^{\infty} \Delta_i$ (here, $S$ is $\Gamma_0\Rightarrow \Delta_0$). Consider the \qciore-structure $\ff A=\langle A, \tripc{\cdot}\rangle$, defined as follows : 
\begin{itemize}
\item[-] $A$ is the set of all free variables occurring in $\Gamma \cup \Delta$ .
\item[-] if $R\in\mm P_n$, then
$$ R^{\ff A}(a_{i_1}, \dots, a_{i_n})=\left\{\begin{array}{ll} \1 & \mbox{ if } R(a_{i_1}, \dots, a_{i_n}) \in \Gamma, \neg R(a_{i_1}, \dots, a_{i_n}) \not\in \Gamma \\
& \\
\um & \mbox{ if } R(a_{i_1}, \dots, a_{i_n}) \in \Gamma, \neg R(a_{i_1}, \dots, a_{i_n}) \in \Gamma \\
& \\
\0 & \mbox{ if } R(a_{i_1}, \dots, a_{i_n}) \not\in \Gamma \end{array} \right.$$
\end{itemize}
and let $s:{\cal V}_{f}\cup {\cal V}_{b}\fun A$ \, be the assignment defined by $s(a)=a$ if $a$ is a free variable, $s(x)$ arbitrary if $x$ is a bound variable. In order to prove that $S$ is not satisfied in $\ff A$ we prove that $s$ satisfies every formula in $\Gamma$ and none in $\Delta$. Let $\alpha$ be a formula in $\Gamma\cup\Delta$. We use induction on the complexity of $\alpha$.\\[2mm]
{\bf Base Case:} 
\begin{itemize}
\item $\alpha$ is $R(a_{i_1}, \dots, a_{i_n})$ for some predicate symbol $R$.
\begin{itemize}
\item If $\alpha\in \Gamma$ then $ R(a_{i_1}, \dots, a_{i_n})\in\Gamma$ and $ R^{\ff A}(a_{i_1}, \dots, a_{i_n})\in \{\1, \um\}$. Therefore $s$ satisfies $\alpha$.
\item If $\alpha\in \Delta$  then $ R(a_{i_1}, \dots, a_{i_n})\in\Delta$ and $ R(a_{i_1}, \dots, a_{i_n})\not\in\Gamma$ then $ R^{\ff A}(s(a_{i_1}), \dots, s(a_{i_n}))=\0$. Therefore $s$ does not satisfy $\alpha$.
\end{itemize}
\item If $\alpha$ is $\neg R(a_{i_1}, \dots, a_{i_n})$ for some predicate symbol $R$.
\begin{itemize}
\item If $\alpha\in \Gamma$ then $\neg R(a_{i_1}, \dots, a_{i_n})\in\Gamma$\\
{\bf Case I:} $R(a_{i_1}, \dots, a_{i_n})\in\Gamma$. Then $ R^{\ff A}(s(a_{i_1}), \dots, s(a_{i_n}))=\um$ and so $s \in \por{R(a_{i_1}, \dots, a_{i_n})}$. That is, $s\in\por{\neg R(a_{i_1}, \dots, a_{i_n})}$ and $s$ satisfies $\neg R(a_{i_1}, \dots, a_{i_n})$. Therefore $s$ satisfies $\alpha$.\\
{\bf Case II:}  $R(a_{i_1}, \dots, a_{i_n})\not\in\Gamma$. Then $ R^{\ff A}(s(a_{i_1}), \dots, s(a_{i_n}))=\0$ and then $s \in \menos{R(a_{i_1}, \dots, a_{i_n})}$. So $s \in \mas{\neg R(a_{i_1}, \dots, a_{i_n})}$ and $s$ satisfies $\neg R(a_{i_1}, \dots, a_{i_n})$. Therefore $s$ satisfies $\alpha$.
\item If $\alpha\in \Delta$. Then $ \neg R(a_{i_1}, \dots, a_{i_n})\in\Delta$ and $ \neg R(a_{i_1}, \dots, a_{i_n})\not\in\Gamma$.\\
{\bf Case I:} $R(a_{i_1}, \dots, a_{i_n})\in\Gamma$. Then $ R^{\ff A}(s(a_{i_1}), \dots, s(a_{i_n}))=\1$ and then $s \in \mas{R(a_{i_1}, \dots, a_{i_n})}$. So $s \in \menos{\neg R(a_{i_1}, \dots, a_{i_n})}$  and $s$ does not satisfy $\neg R(a_{i_1}, \dots, a_{i_n})$. Therefore $s$ does not satisfy $\alpha$.\\
{\bf Case II:} The case $R(a_{i_1}, \dots, a_{i_n})\not\in\Gamma$ is discarded in virtue of reduction step 3 of Definition~\ref{defi6.6}.
\end{itemize}
\end{itemize}
{\bf Inductive step.} Let $\alpha \in \Gamma\cup\Delta$ with complexity of $c(\alpha)=k> 1$. We analyze a few cases, the rest are left to the reader.
\begin{itemize}
\item $\alpha$ is $\alpha_{1}\wedge\alpha_{2}$
\begin{itemize}
\item If $\alpha\in\Gamma$. Then $\alpha_{1}\wedge\alpha_{2}\in\Gamma$. Let $i<\omega$ be the least natural such that $\alpha_{1}\wedge\alpha_{2}\in\Gamma_{i}$. By the reduction step 5 of Definition \ref{defi6.6}, there exists $j>i$ such that $\alpha_{1} , \alpha_{2} \in \Gamma_{j}$ and then $\alpha_{1} , \alpha_{2} \in \Gamma$. By the inductive hypothesis, $s$ satisfies both $\alpha_{1}$ and $\alpha_{2}$. That is, $s \in \mas{\alpha_{1}} \cup \por{\alpha_{1}}$ and $s \in \mas{\alpha_{2}} \cup \por{\alpha_{2}}$\\
{\bf Case I:} $ s \in \mas {\alpha_{1}}$ and $ s \in \mas {\alpha_{2}}$. Then $s \in \mas {\alpha_{1}} \wedge \mas {\alpha_{2}}$ and so $s \in \mas{\alpha_{1}\wedge\alpha_{2}}$. That is $s$ satisfies $\alpha_{1}\wedge\alpha_{2}$. \\
{\bf Case II:} $s \in \mas {\alpha_{1}}$ and $ s \in \por {\alpha_{2}}$. Then $s \in \mas {\alpha_{1}} \wedge \por {\alpha_{2}}$ and so $s \in \mas{\alpha_{1}\wedge\alpha_{2}}$. That is $s$ satisfies $\alpha_{1}\wedge\alpha_{2}$.\\
{\bf Case III:} $s \in \por {\alpha_{1}}$ and $ s \in \mas {\alpha_{2}}$. That is $s \in \por {\alpha_{1}} \wedge \mas {\alpha_{2}}$ and then $s \in \mas{\alpha_{1}\wedge\alpha_{2}}$. Then $s$ satisfies $\alpha_{1}\wedge\alpha_{2}$.\\
{\bf Case IV:} $s \in \por {\alpha_{1}}$ and $s \in \por {\alpha_{2}}$. Then $s \in \por {\alpha_{1}} \wedge \por {\alpha_{2}}$. So $s \in \por{\alpha_{1}\wedge\alpha_{2}}$ and $s$ satisfies $\alpha_{1}\wedge\alpha_{2}$.

\item If $\alpha\in\Delta$. Then $\alpha_{1}\wedge\alpha_{2}\in\Delta$. Let $i<\omega$ be the least natural such that $\alpha_{1}\wedge\alpha_{2}\in\Delta_{i}$. By the reduction step 6 of Definition \ref{defi6.6},  there exists $j >i$ such that either $\alpha_{1}\in\Delta_{j}$ or $\alpha_{2}\in\Delta_{j}$. Suppose that $\alpha_1\in\Delta_{j}$ then $\alpha_{1}\in\Delta$. By the inductive hypothesis, $s$ does not satisfy $\alpha_{1}$ that is $s \in\menos{\alpha_{1}}$ i.e. $s\in\menos{\alpha_{1}}\cup\menos{\alpha_{2}}$. Then $s\in\menos{\alpha_{1}\wedge\alpha_{2}}$ that is $s$ does not satisfy $\alpha_{1}\wedge\alpha_{2}$.
\end{itemize}

\item $\alpha$ is $\neg{\circ}\alpha_1$.

\begin{itemize}
\item If $\alpha\in\Gamma$. Then $\neg{\circ}\alpha_{1}\in\Gamma$. Let $i<\omega$ be the least natural such that $\neg{\circ}\alpha_{1}\in\Gamma_{i}$. By the reduction step 4 of Definition \ref{defi6.6}, there exists $j>i$ such that $\alpha_1 , \neg\alpha_1 \in \Gamma_{j}$ and then $\alpha_1 , \neg\alpha_1 \in \Gamma$. By the inductive hypothesis, $s$ satisfies both $\alpha_1$ and $\neg\alpha_{1}$. That is, $s \in \mas{\alpha_{1}} \cup \por{\alpha_{1}}$ and $s \in \mas{\neg\alpha_{1}} \cup \por{\neg\alpha_{1}}$; that is $s \in\por{\alpha_{1}}$. Then $s \in\menos{{\circ}\alpha_{1}}$ and so $s \in\mas{\neg{\circ}\alpha_{1}}$. Therefore $s$ satisfies $\neg{\circ}\alpha_{1}$. 
\item If $\alpha\in\Delta$. Then $\neg{\circ}\alpha_1\in\Delta$. Let $i<\omega$ be the least natural such that $\neg{\circ}\alpha_1\in\Delta_{i}$. By the reduction step 3 of Definition \ref{defi6.6},  there exists $j >i$ such that ${\circ}\alpha_1\in\Gamma_{j}$ and then ${\circ}\alpha_1\in\Gamma$. By the induction hypothesis, $s$ satisfies ${\circ}\alpha_1$. Then $s \in\mas{{\circ}\alpha_1}\cup \por{{\circ}\alpha_1}$ and, since $\por{{\circ}\alpha_1}=\emptyset$, $s \in\mas{{\circ}\alpha_1}$. Then,  $s \in\menos{\neg{\circ}\alpha_1}$ that is $s$ does not satisfy $\neg{\circ}\alpha_1$.
\end{itemize}

\item If $\alpha$ is $\circ\exists x \alpha_{1}(x)$
\begin{itemize}
\item If $\alpha\in\Gamma$. Then $\circ\exists x \alpha_{1}(x) \in\Gamma$. Let $i<\omega$ be the least natural such that $\circ\exists x \alpha_{1}(x)\in\Gamma_{i}$. By the reduction step 25 of Definition \ref{defi6.6}, there exists $j>i$ and  a free variable symbol $a$ (which are not available at the stage of the formation of $\Gamma_i\Rightarrow \Delta_i$) such that ${\circ}\alpha_{1}(a)\in\Gamma_j$. Then ${\circ}\alpha_{1}(a)\in\Gamma$. By the inductive hypothesis, $s$ satisfies ${\circ}\alpha_{1}(a)$. Then $s \in\mas{{\circ}\alpha_{1} (a)}\cup\por{{\circ}\alpha_{1} (a)}$. But $\por{{\circ}\alpha_{1} (a)}=\emptyset$ and so $s\in \mas{{\circ}\alpha_{1} (a)}$. Then $s\in\mas{\alpha_{1}(a)}\cup\menos{\alpha_{1}(a)}$. Let us suppose that s does not satisfy ${\circ}\exists x\alpha_{1}(x)$ then $s \in\menos{{\circ}\exists x\alpha_{1}(x)}$ and so $s \in\por{\exists x\alpha_{1}(x)}$. Then $ s\in \widehat{\pt a}(\por{\alpha_{1}(a)})$ and so $s_{a}^b\in\por{\alpha_{1}(a)}$ for every free variable $b$, in particular $s\in \por{\alpha_{1}(a)}$ which is a contradiction since $s\in \mas{\alpha_{1}(a)}\cup\menos{\alpha_{1}(a)}$. Therefore $s$ satisfies ${\circ}\exists x \alpha_{1}(x)$.

\item If $\alpha\in\Delta$. Then ${\circ}\exists x \alpha_{1}(x)\in\Delta$. Let $i<\omega$ be the least natural such that ${\circ}\exists x \alpha_{1}(x)\in\Delta_{i}$. By the reduction step 26 of Definition \ref{defi6.6},  there exists $j >i$ such that ${\circ}\alpha_{1}(a)\in\Delta_{j}$ and then ${\circ}\alpha_{1}(a)\in\Delta$ for any free variable $a$ in $A$. By the inductive hypothesis, $s$ does not satisfy ${\circ}\alpha_{1}(a)$ that is $s\in \menos{{\circ}\alpha_{1}(a)}=\por{\alpha_{1}(a)}$ for any free variable $a$. Then $ s\in \widehat{\pt a}(\por{\alpha_{1}(a)})$ and so $s\in \por{\exists x\alpha_{1}(x)}$ and  $s\in \menos{{\circ}\exists x\alpha_{1}(x)}$. Therefore, $s$ does not satisfies ${\circ}\exists x\alpha_{1}(x)$.
\end{itemize}

\item If $\alpha$ is $\circ\forall x \alpha_{1}(x)$
\begin{itemize}
\item If $\alpha\in\Gamma$.  Then $\circ\forall x \alpha_{1}(x) \in\Gamma$. Let $i<\omega$ be the least natural such that $\circ\forall x \alpha_{1}(x)\in\Gamma_{i}$. By the reduction step 23 of Definition \ref{defi6.6}, there exists $j>i$  such that $\circ\alpha_{1}(a)\in\Gamma_{j}$ for any free variable $a$ in $A$.  Then ${\circ}\alpha_{1}(a)\in\Gamma$. By the inductive hypothesis, $s$ satisfies ${\circ}\alpha_{1}(a)$ and so $s \in\mas{{\circ}\alpha_{1} (a)}\cup\por{{\circ}\alpha_{1} (a)}$. But $\por{{\circ}\alpha_{1} (a)}=\emptyset$ and so $s\in \mas{{\circ}\alpha_{1} (a)}$. Then $s\in\mas{\alpha_{1}(a)}\cup\menos{\alpha_{1}(a)}$.
\begin{itemize}
\item If $s\in\mas{\alpha_{1}(a)}$. Then $s\not\in \por{\alpha_1(a)}$ and then $s\not\in \widehat{\pt a}\por{\alpha_1(a)}$. That is $s\not\in \por{\forall x \alpha_1(x)}$ and so $s\in \mas{\forall x \alpha_1(x)}\cup\menos{\forall x \alpha_1(x)}$. Therefore, $s\in \mas{{\circ}\forall x \alpha_1(x)}$ and $s$ satisfies $\circ\forall x \alpha_{1}(x)$.

\item If $s\in\menos{\alpha_{1}(a)}$. Then $s=s_a^a\in\menos{\alpha_{1}(a)}$. That is $s\in \widehat{\ex a}(\menos{\alpha_{1}(a)})$ and so $s\in\menos{\forall x\alpha_{1}(x)}$. Then $s\in\mas{{\circ}\forall x\alpha_{1}(x)}$ and so $s$ satisfies ${\circ}\forall x\alpha_{1}(x)$.  
\end{itemize}

\item If $\alpha\in\Delta$. Then ${\circ}\forall x \alpha_{1}(x)\in\Delta$. Let $i<\omega$ be the least natural such that ${\circ}\forall x \alpha_{1}(x)\in\Delta_{i}$. By the reduction step 24 of Definition \ref{defi6.6},  there exists $j >i$ such that ${\circ}\alpha_{1}(a)\in\Delta_{j}$ for any free variable $a$ in $A$. Then ${\circ}\alpha_{1}(a)\in\Delta$, for any free variable $a$ in $A$. By the inductive hypothesis, $s$ does not satisfy ${\circ}\alpha_{1}(a)$ and therefore  $s\in\menos{{\circ}\alpha_{1}(a)}=\por{\alpha_{1}(a)}$ for any free variable $a$.  Then $ s\in \widehat{\pt a}(\por{\alpha_{1}(a)})$ and so $s\in \por{\forall x\alpha_{1}(x)}$ i.e.  $s\in \menos{{\circ}\forall x\alpha_{1}(x)}$. That is, $s$ does not satisfy ${\circ}\forall x \alpha_{1}(x)$.
\end{itemize}
\end{itemize}
The remaining cases are analyzed analogously.
\end{dem}

\begin{cor}(Completeness Theorem) A sequent $S$ is provable in \gqciore\ if and only if it is valid. 
\end{cor}
\begin{dem} The ``if'' part is consequence of Theorem \ref{soundQ}. For the ``only if'' part, since $S$ is valid the reduction tree $\T(S)$ of $S$ is finite and then we can easily construct a (cut-free) proof (in \gqciore) for $S$.  
\end{dem}

\begin{cor}(Cut elimination property) \gqciore\ admits cut elimination.
\end{cor}

\section{Concluding remarks}\label{s7}

For the purposes of building {\bf LFI}-based theorem provers for real-life applications, it is important to develop  proof theory of the  first-order versions of such {\bf LFI}'s. We think that our work makes a first step in that direction. In this opportunity, we undertook the study of the propositional logic \ciore\ and its first-order version \qciore\ from a proof-theoretic point of view. This logic  was developed as a suitable tool for dealing with inconsistent databases from the point of view of {\em Logics of Formal Inconsistency} ({\bf LFI}s). As it was mentioned, \ciore\ has singular properties that make it an interesting subject of study. In first place, we present a syntactic version of \ciore\ by means of suitable sequent system. Then, we provide a semantic proof of the fact that such system enjoys the cut-elimination property and, as an application, we show some properties of \ciore. Later, we extend the above-mention sequent system to first-order languages providing the sequent rules which govern the behavior of the quantifiers as well as their interaction with the consistency operator. Finally, we prove the completeness and cut-elimination theorem using the well-known Sh\"utte's  method.

We leave for a future work the extension of the first-order system above-mentioned for \qciore\ with equality. Besides, it would be interesting to find applications of the cut-elimination theorem in order to provide syntactic proofs for important results such as the Craig's interpolation theorem, Robinson's theorem and Beth's definability theorem.

\end{document}